\documentclass[11pt,reqno]{article}
\usepackage{amssymb,amsfonts,amsmath,color,url,bm}

\usepackage[english]{babel}

\renewcommand{\le}{\leqslant}
\renewcommand{\leq}{\leqslant}
\renewcommand{\ge}{\geqslant}

\usepackage{latexsym}
\usepackage{amsfonts}
\usepackage{amssymb}
\usepackage{mathrsfs}
\usepackage[all]{xy}
\usepackage{bbm}

\newtheorem{theorem}{Theorem}[section]
\newtheorem{corollary}[theorem]{Corollary}
\newtheorem{lemma}[theorem]{Lemma}
\newtheorem{proposition}[theorem]{Proposition}
\newtheorem{kg}{Theorem $\mathbf{A}$}

\newtheorem{kgg}{Theorem $\mathbf{B}$}

\numberwithin{equation}{section}

\newcommand{\RR}{\mathbb{R}}
\newcommand{\NN}{\mathbb{N}}
\newcommand{\QQ}{\mathbb{Q}}
\newcommand{\ZZ}{\mathbb{Z}}
\newcommand{\mt}{\mapsto}
\newcommand{\cI}{\mathcal{I}}
\newcommand{\cT}{\mathcal{T}}
\newcommand{\cS}{\mathcal{S}}
\newcommand{\cR}{\mathcal{R}}

\newcommand{\cP}{\mathscr{P}}
\newcommand{\sV}{\mathscr{V}}
\newcommand{\cQ}{\mathcal{Q}}
\newcommand{\cV}{\mathcal{V}}
\newcommand{\cA}{\mathbf{A}}
\newcommand{\cB}{\mathbf{B}}
\newcommand{\cU}{\mathcal{U}}

\newcommand{\Bad}{\mathbf{Bad}}

\newcommand{\suc}{\mathrm{suc}}
\newcommand{\CCC}{\mathcal{C}}
\renewcommand{\L}{\mathrm{L}}

\begin{document}

\title{Badly approximable points on  planar curves \\ and  winning}

\author{
Jinpeng An\footnote{Research supported by NSFC grant 11322101 }
\and
Victor Beresnevich\footnote{Research supported by EPSRC grant EP/J018260/1}
\and
 Sanju Velani\footnote{Research supported by EPSRC grant EP/J018260/1}
}

\date{}

\maketitle

\begin{abstract}
For any $i,j  >  0$ with $i+j =1$, let  $\Bad(i,j)$ denote
the set of  points $(x,y) \in \RR^2$ such that $ \max \{
\|qx\|^{1/i}, \, \|qy\|^{1/j} \} > c/q $ for some positive constant $c = c(x,y)$ and all $ q  \in \NN $.  We show that  $\Bad(i,j) \cap \CCC$ is winning in the sense of Schmidt games for
a large class of planar curves $\CCC$, namely, everywhere non-degenerate planar curves and straight lines satisfying a natural Diophantine condition. This
strengthens recent results solving a problem  of
Davenport from the sixties. In short, within the context of Davenport's problem, the winning statement
is best possible. Furthermore, we obtain the inhomogeneous generalizations of the winning results for planar curves and lines and also show that the inhomogeneous form of $\Bad(i,j)$ is winning for two dimensional Schmidt games.
\end{abstract}

\maketitle

{\small

\noindent\emph{Key words and phrases}:
Inhomogeneous Diophantine approximation, non-degenerate curves, Schmidt games

\medskip

\noindent\emph{AMS Subject classification}: 11J83, 11J13, 11K60
}

\section{Introduction}

A real number $x$ is said to be {\em badly approximable} if there
exists a positive constant $c(x)$ such that
\begin{equation*}
 \|qx\| \ > \ c(x) \ q^{-1}  \quad \forall \  q \in \NN   \ .
\end{equation*}
Here and throughout $ \| \cdot  \| $ denotes the distance of a real
number to the nearest integer.  It is well known that the set  $\Bad $
of badly approximable numbers is of Lebesgue measure zero but of
maximal Hausdorff dimension; i.e. $ \dim \Bad = 1 $. In higher
dimensions  there are various natural generalizations of $\Bad$.
Restricting our attention to the plane $\RR^2$, given a pair of real
numbers $i$ and $j$ such that
\begin{equation}\label{neq1}
0 <  i,j  <  1   \quad {\rm and \  } \quad  i+j=1  \, ,
\end{equation}
a point $(x,y) \in \RR^2$ is said to be {\em $(i,j)$-badly
approximable} if there exists a positive constant $c(x,y)$ such that
\begin{equation*}
  \max \{ \; \|qx\|^{\frac1i} \, , \ \|qy\|^{\frac1j} \,  \} \ > \
  c(x,y) \ q^{-1} \quad \forall \  q \in \NN   \ .
\end{equation*}
Denote by   $\Bad(i,j)$  the set of $(i,j)$-badly approximable
points in $\RR^2$.  In the case $i=j=1/2$, the set under consideration is the
standard set of simultaneously badly approximable points.
It easily follows from classical results in the theory of metric
Diophantine approximation   that $\Bad(i,j)$ is of (two-dimensional)
Lebesgue measure zero.   Regarding dimension,  it was shown by Schmidt  \cite{Schmidt-1966}   in the vintage year of 1966 that $\dim \Bad(\frac12,\frac12)=2$.   In fact, Schmidt  proved the significantly stronger statement that $ \Bad(\frac12,\frac12)$ is winning in the sense of his now famous  $(\alpha,\beta)$-games  -- see \S\ref{sgrt}. Almost fourty years later it was proved
in \cite{Pollington-Velani-02:MR1911218} that $\dim \Bad(i,j)=2$ and just recently the first author in \cite{Jinpeng2} has shown that  $ \Bad(i,j) $  is in fact winning.  The latter implies that any countable intersection of $  \Bad(i,j)  $  sets is of full dimension and thus provides a clean and direct proof of Schmidt's Conjecture -- see also \cite{Jinpeng, Badziahin-Pollington-Velani-Schmidt}.

Now let  $\CCC $  be  a planar curve.  Without loss of
generality, we assume that $\CCC $ is  given as a graph
$$
\CCC_f:= \{ (x,f(x))  : x \in I \}
$$
for some function $f$ defined on an interval $I \subset \RR$. Throughout we will assume that $f\in C^{(2)}(I)$, a condition that conveniently allows us to define the curvature. Motivated by a problem of Davenport \cite[p.52]{Davenport-64:MR0166154} from the sixties, the following statement regarding the intersection of $  \Bad(i,j)$ sets with any curve $\CCC$ that is not a straight line segment has recently been established \cite{Badziahin-Velani-DAV, Beresnevich-Dav}.

\begin{kg}\label{thmcountable}
Let  $(i_t,j_t)$ be a countable number of pairs of real numbers
satisfying~\eqref{neq1}  and  suppose that
\begin{equation} \label{liminfassump}
\liminf_{t \to \infty }  \min\{ i_t,j_t  \}  > 0   \ .
\end{equation}
Let
$\CCC :=
\CCC_f$ be  a $C^{(2)}$ planar curve that is not a straight line segment. Then
$$
\dim \Big(\bigcap_{t=1}^{\infty} \Bad(i_t,j_t)   \cap \CCC
\Big) = 1 \ .
$$
\end{kg}

\noindent The theorem implies that there are continuum many points on the parabola $ \cV_2 := \{(x,x^2): x \in \RR\}$ that are simultaneously badly approximable in the $ \Bad(\frac12,\frac12)$ sense and thus provides a solution to the specific problem raised by Davenport in \cite{Davenport-64:MR0166154}. It is worth mentioning that a consequence of \cite[Theorem 1]{B-V-Note} is that the set  $ \Bad(i,j)   \cap \CCC$ is  of zero (induced) Lebesgue measure on $\CCC$. Thus, the fact that it is a set of  full dimension is not trivial.

The condition imposed on $\CCC$ is natural since the statement is not true for all lines.  Indeed,  let
$\L_{a}$ denote the vertical line parallel to the $y$-axis
passing through the point $(a,0)$ in the  $(x,y)$-plane.  Then,
it is  easily  verified, see
\cite[\S1.3]{Badziahin-Pollington-Velani-Schmidt} for the details,
that
$
  \Bad(i,j)  \cap  \L_{a} = \emptyset
$
for any $a \in \RR $ satisfying $ \liminf_{q \to \infty}
q^{1/i}  \|q a\| = 0 \, . $     On the
other hand, if the $\liminf $ is strictly positive then
$
\dim (  \Bad(i,j) \cap \L_{a})  = 1
$.
This  is much harder to prove and is at the heart of the original proof
of Schmidt's Conjecture  established in
\cite{Badziahin-Pollington-Velani-Schmidt}.  Subsequently, it was shown in \cite{Jinpeng} that $ \Bad(i,j) \cap \L_{a} $  is winning.  The following non-trivial  extension of the full dimensional result to non-vertical lines has recently been established  in  \cite{Badziahin-Velani-DAV}.

\begin{kgg}\label{thmcountable2}
Let  $(i_t,j_t)$ be a countable number of pairs of real numbers
satisfying~\eqref{neq1} and~\eqref{liminfassump}. Given $a,
b \in \RR$, let $\L_{a,b} $ denote the  line defined
by the equation $y = a x +b $.   Suppose there exists
$\epsilon>0$ such that
\begin{equation} \label{diocond}
\liminf_{q \to \infty}   q^{\frac{1}{\sigma}-\epsilon}  \|q a
\|
> 0 \, \qquad{\rm where } \ \  \sigma := \sup\{\min\{i_t,j_t\} : t
\in \NN \}.
\end{equation}
Then
$$
\dim \Big(\bigcap_{t=1}^{\infty} \Bad(i_t,j_t)   \cap \L_{a,b} \Big) = 1 \ .
$$
\end{kgg}

Both  Theorem A and Theorem B should be true without the $\liminf$ condition
\eqref{liminfassump}. Indeed, as pointed out in Remark 2 of \cite[\S1.2]{Badziahin-Velani-DAV}, it is very tempting and
not at all outrageous to assert that {\em $ \Bad(i,j)   \cap \CCC $
is winning} at least on the part   of the curve that is genuinely curved.  If true it would imply Theorem A without
assuming \eqref{liminfassump}.   In short, this is precisely what we show in this paper.  We also obtain a winning statement for non-vertical lines that not only implies  Theorem B  without  assuming \eqref{liminfassump} but replaces condition (\ref{diocond}) by a weaker and essentially optimal Diophantine condition.
Furthermore, by making use of a simple idea introduced in \cite{3probs} that provides a natural mechanism for   generalizing    homogenous badly approximable statements to the inhomogeneous setting,  we establish the inhomogeneous  generalization of our winning results.   The same idea is also exploited  to prove that  inhomogeneous $\Bad(i,j)$ is winning.

\subsection{Inhomogeneous $\Bad(i,j)$ and our results  \label{the results}}

For
$\bm\theta=(\gamma,\delta)\in\RR^2$, let $\Bad_{\bm\theta}(i,j)$ denote the set of points  $ (x,y) \in \RR^2 $ such that
\begin{equation*}
  \max \{ \; \|qx - \gamma \|^{\frac1i} \; , \ \|qy - \delta \|^{\frac1j} \,  \} \ > \
  c(x,y) \ q^{-1} \quad \forall \  q \in \NN   \ .
\end{equation*}
It is straight forward to deduce that $\Bad_{\bm\theta}(i,j)$ is of measure zero from the inhomogeneous version of Khintchine's theorem with varying  approximating functions in each co-ordinate. Surprisingly, the fact that $\dim \Bad_{\bm\theta}(i,j) =2  $  is very much a recent development  -- see \cite{ET} for the symmetric $i=j=1/2$ case and \cite{3probs}  for the general case.

One of the main goals of this paper is to prove the   following full dimension statement which not only implies the inhomogeneous analogue of Theorem A but totally  removes the   $\liminf$ condition \eqref{liminfassump}.

\begin{theorem}\label{thm1}
Let  $(i,j)$ be a pair of real numbers
satisfying~\eqref{neq1}  and let
$\CCC :=
\CCC_f$ be  a planar curve such that  $f \in C^{(2)}(I) $ and  that  $ f''(x) \neq  0 $ for all $ x \in I $.
Then, for any $\bm\theta=(\gamma,\delta)\in\RR^2$ we have that  $ \Bad_{\bm\theta}(i,j) \cap \CCC $ is a  winning subset of $\CCC$.
\end{theorem}

\noindent{\em Remark 1.}
The condition $f''(x)\neq0$ is often referred to as non-degeneracy at $x$.
Note that if $f \in C^{(2)}(I) $ and  $ f''(x) \neq  0 $ for some point $ x \in I $, then, by continuity, there exists an interval $I^* \subset I$ such that  $ f''(x) \neq  0 $ for all $ x \in I^* $. In other words, if the curve $\CCC_f$ is not a straight line segment then it is always possible to find an arc of $\CCC_f$ that is non-degenerate everywhere.  Thus,  Theorem \ref{thm1} with $I = I^*$  and $\bm\theta =(0,0)$  implies Theorem A without
assuming \eqref{liminfassump}.  Of course this makes use of the well know  fact that any winning set is of full dimension and  that  any countable intersection of winning sets is again winning.

\vspace*{2ex}

\noindent{\em Remark 2.} Given a subset $X \subset \RR^2$, let $\pi(X)$ denote the projection of $X$ onto the
$x$-axis.  Regarding Theorem \ref{thm1},  we actually prove that {\em $\pi (\Bad_{\bm\theta}(i,j) \cap \CCC)   $ is a $1/2$-winning subset of~$I$}. In other words, for the projected set we prove $\alpha$-winning with the best possible winning constant; i.e. $\alpha=1/2$.
Now, since the projection  map $\pi$ is bi-Lipschitz on $\CCC$ and the image of a winning set under  a  bi-Lipschitz map is
again winning  \cite[Proposition 5.3]{DaniDodVic},  it trivially follows that
$\Bad_{\bm\theta}(i,j) \cap \CCC $ is an $ \alpha_0$-winning subset of $\CCC$ for some
$ \alpha_0 \in (0,1/2]$. The actual value of $ \alpha_0 $ is dependent on the Lipschitz constant $\kappa$  associated with $(\pi|_\CCC)^{-1}$.  In fact,  if we use the maximum norm on $\RR^2$ it is  possible to obtain the winning statement  with $ \alpha_0=1/2$.  Essentially, if $\kappa > 1 $ we  consider the projection of $ \Bad_{\bm\theta}(i,j) \cap \CCC $ onto the $y$-axis rather then the  $x$-axis.

\vspace*{2ex}

For straight lines we prove the following  counterpart  statement.

\begin{theorem}\label{thm2}
Let  $(i,j)$ be a pair of real numbers
satisfying~\eqref{neq1} and given $a,b\in \RR$ with $a \neq 0$, let $\L_{a,b } $ denote the  line defined
by the equation $y = a x +b $.   Suppose there exists
$\epsilon>0$ such that
\begin{equation} \label{diocondbetter}
\liminf_{q \to \infty}   q^{\frac{1}{\sigma}-\epsilon}  \max \{\|q a
\|, \|q b
\|  \}
> 0 \, \qquad{\rm where } \ \  \sigma := \min\{i,j\} \, .
\end{equation}
Then, for any $\bm\theta=(\gamma,\delta)\in\RR^2$ we have that   $ \Bad_{\bm\theta} (i,j) \cap \L_{a,b }$ is a winning subset of $\L_{a,b }$.  Moreover, if $ a \in \QQ$ the statement is true with $ \epsilon=0$ in  (\ref{diocondbetter}).
\end{theorem}

\noindent{\em Remark 3.} As  with curves,  the theorem is actually deduced on showing that the projected set
$
\pi(\Bad_{\bm\theta} (i,j) \cap \L_{a,b }) $ is an $1/2$-winning subset of $\RR$.

\vspace*{2ex}

\noindent{\em Remark 4.} The  Diophantine condition (\ref{diocondbetter}) imposed in the theorem  is essentially optimal since
\begin{equation}  \label{lop} \Bad(i,j) \cap \L_{a,b }   = \emptyset   \quad   {\rm  if }   \quad  \liminf_{q \to \infty}   q^{\frac{1}{\sigma}}  \max \{\|q a
\|, \|q b
\|  \} = 0. \end{equation}
To see that this is the case,  assume  for the moment that $\Bad(i,j) \cap \L_{a,b }$ is nonempty.  Then there exists some point  $x\in\RR$ such that $(x,ax+b)\in \Bad(i,j)$.  In terms of the equivalent dual form representation of $\Bad(i,j)$ -- see \cite[\S1.3]{Badziahin-Pollington-Velani-Schmidt}, this means that there exists a constant  $c(x)>0$ such  that
\begin{equation}  \label{op} \max\{|A|^{\frac{1}{i}},|B|^{\frac{1}{j}}\}  \, |Ax+B(ax+b)+C| \ \ge \ c(x)  \, \end{equation}
for all $A,B,C\in\ZZ$ with $(A,B)\ne(0,0)$.    Now, for any given  $B\in\NN$ we choose $A,C\in\ZZ$  such that $|Ba+A|=\|Ba\|$ and $|Bb+C|=\|Bb\|$. Then
\begin{eqnarray*}
|Ax+B(ax+b)+C| & =  &  |(Ba+A)x+(Bb+C)|  \ \le \ \|Ba\|  \,  |x|  \, +   \, \|Bb\| \\[1ex]
& \le   & (1+|x|)\max\{\|Ba\|,\|Bb\|\}
\end{eqnarray*}
and
\begin{eqnarray*}
\max\{|A|^{\frac{1}{i}},|B|^{\frac{1}{j}}\} & \le  &  \max\{(|Ba|+1)^{\frac{1}{i}},|B|^{\frac{1}{j}}\}  \\[1ex] &  \le  &  \max\{(1+|a|)^{\frac{1}{i}}|B|^{\frac{1}{i}},|B|^{\frac{1}{j}}\}   \ \le  \  (1+|a|)^{\frac{1}{i}} |B|^{\frac{1}{\sigma}}.
\end{eqnarray*}
Thus, on combining these estimates with  (\ref{op}), it follows that
$$|B|^{\frac{1}{\sigma}}\max\{\|Ba\|,\|Bb\|\} \ \ge \ \frac{c(x)}{(1+|a|)^{\frac{1}{i}}(1+|x|)}   \quad \forall \ B \in \NN    $$
and so
$$\liminf_{q \to \infty}   q^{\frac{1}{\sigma}}  \max \{\|q a
\|, \|q b
\|  \} > 0. $$
This establishes (\ref{lop}).    Note than in view  of  (\ref{lop}) and the moreover part of the theorem, the Diophantine condition (\ref{diocondbetter}) with $ \epsilon=0$ is optimal in the case $a$ is rational.

\vspace*{2ex}

\noindent{\em Remark 5.}  The fact that $a = 0 $ is excluded in the statement of the theorem is natural  since as in Remark 4, on making use of the equivalent dual form representation of $\Bad(i,j)$ it is  easily verified that
\begin{equation*}  \label{llop} \Bad(i,j) \cap \L_{0,b }   = \emptyset   \quad   {\rm  if }   \quad  \liminf_{q \to \infty}   q^{\frac{1}{j}}   \|q b
\|  = 0. \end{equation*}   On the
other hand, if the above  $\liminf $ is strictly positive then it was shown in \cite{Jinpeng} that $ \Bad(i,j) \cap \L_{0,b} $  is a winning subset of the horizontal line $\L_{0,b}$.  By making use of the mechanism developed in this paper, it is relatively straightforward to adapt the homogeneous  proof given  in \cite{Jinpeng} to show that if $ \liminf_{q \to \infty}   q^{1/j}   \|q b
\|  > 0$,  {\em  then for any $\bm\theta \in \RR^2$ the set  $ \Bad_{\bm\theta}(i,j) \cap \L_{0,b} $  is a winning subset of the horizontal line $\L_{0,b}$.}

\vspace*{2ex}

\noindent{\em Remark 6.}  Observe that when it comes to intersecting countably many  $(i_t,j_t)$ pairs  the Diophantine condition (\ref{diocondbetter}) imposes  the condition that
\begin{equation} \label{diocondbetterohyes}
\liminf_{q \to \infty}   q^{\frac{1}{\sigma}-\epsilon}  \max \{\|q a
\|, \|q b
\|  \}
> 0 \, \qquad{\rm where } \ \  \sigma := \sup\{\min\{i_t,j_t\} : t
\in \NN \}.
\end{equation}
This is clearly weaker than condition (\ref{diocond}) imposed in Theorem B and moreover in view of Remark 4 it   is essentially optimal.

\medskip

The proofs of Theorem \ref{thm1} and Theorem \ref{thm2} rely on first establishing the homogeneous cases and then making use of a  natural  mechanism that we develop  for   generalizing  homogenous winning statements to the inhomogeneous setup.  This mechanism is further exploited to prove that  inhomogeneous $\Bad(i,j)$ is winning.

\begin{theorem}\label{2DT:main}
Let  $(i,j)$ be a pair of real numbers satisfying~\eqref{neq1}. Then, for any $\bm\theta=(\gamma,\delta)\in\RR^2$ we have that  $\Bad_{\bm\theta}(i,j)$ is a  $(30\sqrt{2})^{-1}$-winning subset of $\RR^2$.
\end{theorem}

\noindent{\em Remark 7.} It is worth pointing out that the winning constant $(30\sqrt{2})^{-1}$ is not optimal. Indeed, the ideas used to prove the above  theorem  and  the argument given in \cite{NS} can be combined to show that $\Bad_{\bm\theta}(i,j)$ is hyperplane absolute winning. This is a stronger version of winning and implies that  $\Bad_{\bm\theta}(i,j)$ is $\alpha$-winning for any $\alpha < 1/2$.

\section{The main strategy }

In this section we  outline  the key steps in establishing the homogeneous case ($\bm\theta =(0,0)$) of Theorem \ref{thm1}. The general inhomogeneous statement is obtained by appropriately adapting  the homogeneous argument and is carried out in \S\ref{inhThm1.2}.  To begin with observe that for any planar curve $\CCC :=
\CCC_f$ and  $\bm\theta \in \RR^2$
$$
\Bad_{\bm\theta}^f(i,j)  \, :=  \, \{x\in I:(x,f(x))\in\Bad_{\bm\theta}(i,j)\}  \; =  \; \pi(\Bad_{\bm\theta}(i,j) \cap \CCC).
$$
Recall, that $\pi : \RR^2 \to \RR $  is the projection map onto the $x$-axis.   Also, for convenience and without loss of generality we will assume that $j \leq i$. Thus, the homogeneous case of Theorem \ref{thm1} is easily deduced from  the following  statement for $\Bad^f(i,j):=\Bad_{(0,0)}^f(i,j)$ -- see Remark 2 above for the justification.

\begin{theorem}\label{T:main}
Let  $(i,j)$ be a pair of real numbers
satisfying   $0<j\le i<1$ and $i+j=1$.   Let $I\subset\RR$ be a compact interval and  $f \in C^{(2)}(I) $ such  that  $ f''(x) \neq  0 $ for all $ x \in I $.   Then
$\Bad^f(i,j)$
is a $1/2$-winning subset of $I$.
\end{theorem}

At this point it is useful to recall the definition of a winning set and the notion of rooted trees -- a key `structural' ingredient  in establishing the above winning statement.

\subsection{Schmidt games and rooted trees \label{sgrt}}

Wolfgang M. Schmidt introduced the games which now bear his name in \cite{Schmidt-1966}.   The simplified account which we are about to present  is sufficient for the purposes of this paper.   Suppose that $0<\alpha<1$ and $0 <  \beta < 1$.   Consider the following game involving the two arch rivals $\cA $yesha  and $\cB $hupen --  often simply referred to as players \textbf{A} and \textbf{B}.  First, \textbf{B}  chooses a  closed ball $\cB_0 \subset \RR^m$. Next,  \textbf{A} chooses a closed ball  $\cA_0$ contained in $\cB_0$ of diameter  $ \alpha \, \rho( \cB_0 ) $ where $\rho(\ .\ )$ denotes the diameter of the ball under consideration.  Then,  \textbf{B}  chooses at will a closed ball $\cB_1$ contained in $\cA_0$ of diameter  $ \beta \, \rho( \cA_0 )$.  Alternating in this manner between the two players, generates  a nested sequence of closed balls  in $\RR^m$:
$$
\cB_0\supset \cA_0\supset \cB_1 \supset \cA_1 \supset\ldots \supset \cB_n\supset \cA_n\supset \ldots
$$
with diameters
$$
\rho(\cB_n) \, = \,   (\alpha \, \beta )^{n}  \,   \rho(\cB_0)  \text{\quad and \quad} \rho(\cA_n) \, =  \,   \alpha \,   \rho(\cB_n)      \, .
$$
A subset  $X$ of $\RR^m$ is said to be {\em $(\alpha,\beta)$-winning} if  $\cA$  can play in such a way that the unique point of intersection
$$
\bigcap_{n=0}^{\infty}   \cB_n   \,  = \, \bigcap_{n=0}^{\infty}   \cA_n
$$
lies in $X$, regardless of how  $\cB$ plays.  The set $X$ is called {\em
$\alpha$-winning} if  it is  $(\alpha,\beta)$-winning for all
$\beta\in (0,1)$. Finally, $X$ is simply called {\em winning} if it
is $\alpha$-winning for some~$\alpha$.  Informally, player $\cB$ tries
to stay away from the `target' set $X$ whilst player $\cA$ tries to land
on $X$. As shown by Schmidt \cite{Schmidt-1966}, the following are two key consequences of winning.
\begin{itemize}
\item If $X \subset \RR^m$ is a winning set, then   $\dim X=m$.
\item The intersection of countably many  $\alpha$-winning sets is $\alpha$-winning.
\end{itemize}
In the setting of Theorem \ref{T:main}, we have  $m=1$ and $X=\Bad^f(i,j)$. Thus $\cA_n$ and $\cB_n$ are compact intervals.  Note that more generally, we can replace $\RR^m $ in the above description of Schmidt games by a $m$-dimensional Riemannian  manifold.  It is this slightly more general form that is implicitly referred to within the context of Theorem~\ref{thm1}.

\medskip
We now turn our attention to rooted trees.  Recall that a \emph{rooted tree} is a connected
graph $\cT$ without cycles and with a distinguished vertex $\tau_0$, called the \emph{root} of $\cT$ . We
identify $\cT$ with the set of its vertices. Any vertex  $\tau \in  \cT$ is connected to $\tau_0$ by a unique
path. The length of the path is called the \emph{height} of $\tau$ . The set of vertices of height $n$ is
called the \emph{$n$'th level} of $\cT$ and is denoted by $\cT_n$. Thus $\cT_0= \{ \tau_0 \} $. Next, given $ \tau, \tau' \in \cT$  we write $\tau\prec\tau'$
 to indicate that the path between $\tau_0$  and $\tau$   passes through $\tau'$ and in this case we call  $\tau$ a  \emph{descendant} of $\tau'$ and  $\tau'$   an \emph{ancestor} of $\tau$ . By definition, every vertex is a
descendant and an ancestor of itself. For  $\cV   \subset \cT$ , we write $\tau\prec \cV$ if $\cV$ contains an ancestor
of $\tau$. If $\tau\prec\tau'$ and the height of $\tau$  is one greater than that of $\tau'$, then $\tau$ is called a \emph{successor}
of $\tau'$  and $\tau'$  is called the \emph{predecessor} of $\tau$. Let $\cT (\tau )$ denote the rooted tree formed by
all descendants of $\tau$ . Thus  the  root of $\cT (\tau ) $ is  $\tau$   and we  denote by $\cT_{\suc}( \tau )  $ the set of
all successors of $\tau$ . More generally, for $\cV \subset \cT $, we let  $\cT_{\suc}(\cV):=
\bigcup_{\tau  \in \cV}  \cT_{\suc}(\tau)$. In this
paper, we use the convention that a subtree of $\cT$  has the same root as  $\cT$.  As a consequence,  $\cT (\tau)$  is not regarded as a subtree of $\cT$ unless $\tau  = \tau_0$.

Let $ N  \in \NN$. We say that a rooted tree is \emph{N-regular} if every vertex has exactly $ N $
successors. Note that an $N$-regular rooted tree is necessarily infinite. The following
statement  appears as Proposition 2.1 in  \cite{Jinpeng}.

\begin{proposition} \label{propJA}
Let $\cT$ be an $N$-regular rooted tree, $\cS \subset  \cT$  be a subtree, and $1 \le m \le N$
be an integer. Suppose that for every $m$-regular subtree $\cR$ of $\cT$,  we have that $\cS \cap \cR $ is infinite. Then $\cS$  contains a  $(N -m + 1)$-regular subtree.
\end{proposition}

\noindent This proposition will be required in establishing Proposition \ref{P:main}  below.  As shown in \S\ref{howzthat}, the latter is very much at the heart of the proof of   Theorem \ref{T:main}.

\subsection{The winning strategy for Theorem \ref{T:main} \label{winstrat}}
Let $\beta\in(0,1)$. We want to prove that $\Bad^f(i,j)$ is $(\frac{1}{2},\beta)$-winning. In the first round of the game, \textbf{B}hupen chooses a closed interval $\cB_0\subset I$.  Now  \textbf{A}yesha chooses the closed interval $\cA_0\subset\cB_0$ with diameter $\rho(\cA_0)=\frac{1}{2} \rho(\cB_0)$ such that $\cA_0$ has the same center as $\cB_0$. Let $\kappa>1$ be sufficiently large so  that for every $x\in I$ we have that
 \begin{equation}\label{sv1} |f'(x)|\le \kappa-1
 \end{equation}
 and
\begin{equation}\label{sv2}
 \kappa^{-1}\le|f''(x)|\le \kappa  \, .
 \end{equation}
Clearly such a $\kappa> 1$ exists by the conditions imposed on $f$ and  $I$.    Let
$$R:=(2\beta^{-1})^5  \qquad {\rm and }  \qquad l:=  \rho(\cA_0)  \, . $$
Without loss of generality, we may assume that
\begin{equation}\label{E:l1}
3\kappa lR^2<1.
\end{equation}
The point is that if this is not the case then \textbf{A}yesha  will choose her intervals $\cA_{n-1}   $  $ (n\ge1)$ arbitrarily until    the interval $\cB_n$ chosen by \textbf{B}hupen  with diameter $\rho(\cB_n)  =  (\beta/2)^n  \rho (\cB_0) $ satisfies $ 6\kappa\rho(\cB_n) R^2<1$.  At this stage  \textbf{A}yesha  chooses the closed interval $\cA_n$  with the same center as $\cB_n$ and  half its diameter. We  now simply relabel $\cA_n$ and $\cB_n$  as  $\cA_0$ and $\cB_0$ respectively.

Choose  $\mu>0$   such that
\begin{equation}\label{E:mu}
10\kappa^2l^{-1}R^{\frac{1}{j}-\frac{j}{6}\mu}\le1
\end{equation}
and define
\begin{equation}\label{vb201}
\lambda_0:=0   \qquad  {\rm and }  \qquad \lambda_k:=\frac{2(1+i)}{j} k+\mu  \quad {\rm for }  \quad  k\ge1.
\end{equation}
In turn,  let
\begin{equation} \label{E:c}
c:=\frac{l^2}{10^3\kappa^5R^{4+\lambda_1}}
\end{equation}
and
\begin{equation} \label{E:cP}
\cP:=\Big\{P=\Big(\frac{p}{q},\frac{r}{q}\Big):\frac{p}{q}\in I,\Big|f\Big(\frac{p}{q}\Big)-\frac{r}{q}\Big|<\frac{\kappa c}{q^{1+j}}\Big\}.
\end{equation}
Note that $R$ is large while $l$ and $c$ are small. Indeed, we have the inequalities
$$
R>32,\qquad l<10^{-3},\qquad c<10^{-10}l
$$
which we will use without further reference. Throughout, when a rational
point in $\RR^2$  is expressed as $ P=(\frac{p}{q},\frac{r}{q})$, we assume that $q > 0 $ and  that the integers $ p, q, r $ are co-prime.   Finally, for  each   $P=(\frac{p}{q},\frac{r}{q})\in \QQ^2$  we associate the interval
\begin{equation} \label{E:Delta}
\Delta(P):=\Big\{x\in I:\Big|x-\frac{p}{q}\Big|<\frac{c}{q^{1+i}}\Big\}.
\end{equation}

The following inclusion is a simple consequence of the manner in which the above quantities  and   objects have been defined.

\begin{lemma}\label{L:Bad_c}  Let $\cA_0$, $  \cP  $ and $ \Delta(P) $ be as above.  Then
$$\cA_0\setminus\bigcup_{P\in\cP}\Delta(P) \ \subset \ \Bad^f(i,j)  \, . $$
\end{lemma}

\noindent{\em Proof. \ }
Let $x\in\cA_0$. Suppose $x\notin\Bad^f(i,j)$. Then there exists $P=(\frac{p}{q},\frac{r}{q})\in\QQ^2$ such that $$\Big|x-\frac{p}{q}\Big|<\frac{c}{q^{1+i}}, \qquad \Big|f(x)-\frac{r}{q}\Big|<\frac{c}{q^{1+j}}.
$$
In view of the fact that
$$\Big|x-\frac{p}{q}\Big|<\frac{c}{q^{1+i}}\le c\le\frac{l}{2},$$
it follows that $\frac{p}{q}\in\cB_0\subset I$. Hence, using the Mean Value Theorem together with \eqref{sv1} we obtain the following estimate:
\begin{eqnarray}
\Big|f\Big(\frac{p}{q}\Big)-\frac{r}{q}\Big| & \le  & \Big|f\Big(\frac{p}{q}\Big)-f(x)\Big|+
\Big|f(x)-\frac{r}{q}\Big| \nonumber \\[2ex]
&\le  & (\kappa-1)\Big| x-\frac{p}{q}\Big|+\frac{c}{q^{1+j}} \nonumber \\[2ex]
& < & \frac{(\kappa-1)c}{q^{1+i}}+\frac{c}{q^{1+j}} \ \le  \ \frac{\kappa c}{q^{1+j}}  \, \nonumber .
\end{eqnarray}
The upshot is that $x\in\Delta(P)$ with $P\in\cP$. This completes the proof of the lemma.
\vspace*{-2ex}

\hfill $ \boxtimes $

\vspace*{2ex}

Now let $\cT$ be an $[R]$-regular rooted tree with root $\tau_0$, where $[\ \cdot \ ]$ denotes the integer part of a real number. We choose and fix
an injective map $\cI$ from $\cT$ to the set of closed subintervals of $\cA_0$ satisfying the following conditions:
\begin{itemize}
  \item For any $n\ge0$ and $\tau\in\cT_n$, $\rho(\cI(\tau))=lR^{-n}$. In particular, $\cI(\tau_0)=\cA_0$.
  \item For $\tau,\tau'\in\cT$, if $\tau\prec\tau'$, then $\cI(\tau)\subset\cI(\tau')$.
  \item For any $\tau\in\cT$, the interiors of the intervals $\{\cI(\tau'):\tau'\in\cT_{\suc}(\tau)\}$ are mutually disjoint, and $\bigcup_{\tau'\in\cT_{\suc}(\tau)}\cI(\tau')$ is connected.
\end{itemize}
Note that for $n\ge1$ and $\tau\in\cT_{n-1}$, any closed subinterval of $\cI(\tau)$ of length $2lR^{-n}$ must contain an $\cI(\tau')$ for some $\tau'\in\cT_{\suc}(\tau)$.  Suppose that $\cP$ is partitioned into a disjoint union
$$\cP=\bigcup_{n=1}^\infty\cP_n.$$
We  inductively define a subtree $\cS$ of $\cT$ as follows.
Let $\cS_0=\{\tau_0\}$. If $\cS_{n-1}$ $(n\ge1)$ is defined, we let
$$\cS_n:=\Big\{\tau\in\cT_{\suc}(\cS_{n-1}):\cI(\tau) \ \cap\bigcup_{P\in\cP_n}\Delta(P)=\emptyset\Big\}.$$
Then $$
\cS:=\bigcup_{n=0}^\infty\cS_n$$ is a subtree of $\cT$ and by construction
\begin{equation}\label{E:tau}
\cI(\tau) \ \subset \ \cA_0\setminus\bigcup_{P\in\cP_n}\Delta(P)   \qquad \forall  \ n\ge1 \quad {\rm and } \quad  \tau\in\cS_n.
\end{equation}
Thus given a partition $\cP_n$ of $\cP$,  the intervals $\{\cI(\tau):  \tau \in \cS_n\}$ serve as possible
candidates when it comes to \textbf{A}yesha to turn to make a move.  Moreover, we are able to choose the partition $\cP_n$  in such a way that $\cS$ has the following key feature.
\begin{proposition}\label{P:main}
There exists a partition $\cP=\bigcup_{n=1}^\infty\cP_n$ such that the tree $\cS$ has an $([R]-10)$-regular subtree.
\end{proposition}

Armed with this proposition  we are able to describe the  winning strategy  that  \textbf{A}yesha will adopt.

\subsubsection{Proof of Theorem \ref{T:main} modulo Proposition \ref{P:main} \label{howzthat}}

Let $\cP=\bigcup_{n=1}^\infty\cP_n$ be a partition such that $\cS$ has an $([R]-10)$-regular subtree, say $\cS'$. We inductively prove that for every $n\ge0$,
\begin{equation}\label{E:state}
\text{\textbf{A}yesha can choose $\cA_{5n}=\cI(\tau_n)$ for some $\tau_n\in\cS'_n$.}
\end{equation}
Since $\cA_0=\cI(\tau_0)$, we trivially have that \eqref{E:state} holds when $n=0$. Assume $n\ge1$ and that  \textbf{A}yesha has chosen $\cA_{5(n-1)}=\cI(\tau_{n-1})$, where $\tau_{n-1}\in\cS'_{n-1}$. We refer to  the intervals $\{\cI(\tau):\tau\in\cT_{\suc}(\tau_{n-1})\setminus\cS'_{\suc}(\tau_{n-1})\}$ as dangerous intervals -- they represent  intervals that  \textbf{A}yesha needs to avoid.
We first prove that
\begin{align}\label{E:substate}
&\text{For $t\in\{0,1,2,3,4\}$, \textbf{A}yesha can play so that $\cA_{5(n-1)+t}$, and hence } \notag\\
&\text{$\cB_{5(n-1)+t+1}$, contains at most $[10\cdot2^{-t}]$ dangerous intervals.}
\end{align}
If $t=0$, there is nothing to prove. Assume $1\le t\le 4$ and \eqref{E:substate} holds if $t$ is replaced by $t-1$. Thus $\cB_{5(n-1)+t}$ contains at most $[10\cdot2^{-t+1}]$ dangerous intervals. Divide $\cB_{5(n-1)+t}$ into two closed subintervals of equal length. Then \textbf{A}yesha can choose $\cA_{5(n-1)+t}$ to be one of the subintervals so that it contains at most
$\Big[\frac{1}{2}[10\cdot2^{-t+1}]\Big]\le[10\cdot2^{-t}]$ dangerous intervals. This proves \eqref{E:substate}.
By letting $t=4$ in \eqref{E:substate}, we see that \textbf{A}yesha can play so that $\cB_{5n}$ contains no dangerous intervals. Since $\cB_{5n}$ has length $2lR^{-n}$, it contains an $\cI(\tau_n)$ for some $\tau_n\in\cT_{\suc}(\tau_{n-1})$. It follows that $\tau_n\in\cS'_n$. So \textbf{A}yesha can choose $\cA_{5n}=\cI(\tau_n)$. This completes the proof of \eqref{E:state}.
In view of \eqref{E:state}, \eqref{E:tau} and Lemma \ref{L:Bad_c}, we have
\begin{eqnarray*}
\bigcap_{n=0}^\infty \cA_n & = &\bigcap_{n=1}^\infty \cA_{5n} \ = \ \bigcap_{n=1}^\infty\cI(\tau_n) \ \subset \ \bigcap_{n=1}^\infty\cA_0\setminus\bigcup_{P\in\cP_n}\Delta(P)\\
& = & \cA_0\setminus\bigcup_{P\in\cP}\Delta(P) \ \subset \ \Bad^f(i,j).
\end{eqnarray*}
This proves the theorem assuming  the truth of Proposition \ref{P:main}.
~\vspace*{2ex}
\hfill $ \boxtimes $

\section{Preliminaries for Proposition \ref{P:main}  \label{premain}  }

The following simple but important lemma was established in \cite{Badziahin-Pollington-Velani-Schmidt}.

\begin{lemma}\label{lem31415}
For any point $P=(\frac{p}{q},\frac{r}{q})\in\QQ^2$ there exist coprime integers $A$, $B$, $C$ with
$(A,B)\ne(0,0)$ such that
$$
Ap+Br+Cq=0,
$$
$$
|A|\le q^i  \,   \quad {\rm and }  \quad  |B|\le q^j  \, .
$$
\end{lemma}

\noindent{\em Proof. \ }
Since the proof is only a few lines we reproduce it here. By Minkowski's theorem for systems of linear forms there is $(A,B,C)\in\ZZ^3\setminus\{\mathbf{0}\}$ such that
$$
|Ap+Br+Cq|<1,\quad |A|\le q^i  \quad \text{and}  \quad  |B|\le q^j  \, .
$$
Since $Ap+Br+Cq$ is an integer it must be zero. If $(A,B)=(0,0)$ then $qC=0$ and, since $q\neq0$ we also have that $C=0$, a contradiction. Hence $(A,B)\not=(0,0)$ and the proof is complete.
\vspace*{-2ex}

\hfill $ \boxtimes $

\vspace*{2ex}

In view of Lemma~\ref{lem31415}, to each point $P=(\frac{p}{q},\frac{r}{q})\in\cP$, we can assign a rational line
\begin{equation}\label{vb210}
L_P:=\{(x,y)\in\RR^2:A_Px+B_Py+C_P=0\}
\end{equation}
passing through $P$ where  $A_P,B_P,C_P\in\ZZ$  are co-prime  with  $(A_P,B_P)\ne(0,0)$ and such that
\begin{equation}\label{vb204}
|A_P|\le q^i  \,   \quad {\rm and }  \quad  |B_P|\le q^j  \, .
\end{equation}
If there is more than one line satisfying the above conditions, we choose any one.
Further, for  $P\in\cP$ we define the function $F_P:I\to\RR$ by
$$F_P(x):=A_Px+B_Pf(x)+C_P$$
and the set
\begin{equation}\label{E:Theta}
\Theta(P):=\Big\{x\in I:|F_P(x)|<\frac{2\kappa c}{q}\Big\}.
\end{equation}
In this section we gather basic information regarding the sets $\Theta(P)$ and associated quantities.

\begin{lemma}\label{L:Theta}
Let $P=(\frac{p}{q},\frac{r}{q})\in\cP$ and $x\in\Delta(P)$. Then $$\Big|F_P\Big(\frac{p}{q}\Big)\Big|<\frac{\kappa c}{q}  \qquad {\rm and } \qquad \Big|F_P(x)-F_P\Big(\frac{p}{q}\Big)\Big|<\frac{\kappa c}{q}.$$
In particular, we have that $\Delta(P)\subset\Theta(P)$.
Furthermore, if $x\in\Theta(P)$ then
$$
\Big|F_P(x)-F_P\Big(\frac{p}{q}\Big)\Big|  \ < \ \frac{3\kappa c}{q}.
$$
\end{lemma}

\noindent{\em Proof. \ }   On using the fact that $ P$ lies on the line  $L_P$ and \eqref{E:cP} we obtain the following
$$\Big|F_P\Big(\frac{p}{q}\Big)\Big| \ = \ \Big|B_P\Big(f\Big(\frac{p}{q}\Big)-\frac{r}{q}\Big)\Big|  \ \le  \  q^j \Big|f\Big(\frac{p}{q}\Big)-\frac{r}{q}\Big| \ < \ \frac{\kappa c}{q}   \, . $$
For the second inequality, by the  Mean Value Theorem there exists a point $\xi \in I$ such that
$$\Big|F_P(x)-F_P\Big(\frac{p}{q}\Big)\Big|   =  |F_P'(\xi)| \Big|x-\frac{p}{q}\Big|  \ <  \  (|A_P|+(\kappa-1)|B_P|) \frac{c}{q^{1+i}} \ \le \ \frac{\kappa c}{q}. $$
 In particular, it follows that for $x\in\Delta(P)$,  we have that
 $$|F_P(x)|  \ \le \ \Big|F_P\Big(\frac{p}{q}\Big)\Big| + \Big|F_P(x)-F_P\Big(\frac{p}{q}\Big)\Big| \ < \ \frac{2\kappa c}{q}
 $$
 and so by definition $ x \in \Theta(P)$.

Finally, if $x\in\Theta(P)$, then
$$\Big|F_P(x)-F_P\Big(\frac{p}{q}\Big)\Big|  \ \le \  |F_P(x)|   +  \Big|F_P\Big(\frac{p}{q}\Big)\Big|   \ <  \ \frac{3\kappa c}{q}.\vspace*{-2ex}$$
\hfill $ \boxtimes $

\vspace*{2ex}

Our next goal is to describe the structure of $\Theta(P)$, in particular, to estimate its size. To this end, we introduce the following quantities.
Let
\begin{equation}\label{vb2003}
E_P:=F'_P\Big(\frac{p}{q}\Big)=A_P+B_Pf'\Big(\frac{p}{q}\Big)
\end{equation}
and let
\begin{equation}\label{E:cP*}
\cP^*:=\Big\{P=\Big(\frac{p}{q},\frac{r}{q}\Big)\in\cP:qE_P^2<9\kappa^2c|B_P|\Big\} \,.
\end{equation}
Then we define the height of $P$ by
\begin{equation}\label{vb205}
H(P) \ := \ \max\{q|E_P|,3\kappa\sqrt{cq|B_P|}\} \ = \
\begin{cases}
3\kappa\sqrt{cq|B_P|}  & \ {\rm if }  \quad P\in\cP^* \\[1ex]
q|E_P|   & \ {\rm if }  \quad  P\in\cP\setminus\cP^*  \, .
\end{cases}
\end{equation}

\vspace*{2ex}

\begin{lemma}\label{L:bound}
Let $P=(\frac{p}{q},\frac{r}{q})\in\cP$. Then
\begin{eqnarray}
|E_P| & \le & \kappa q^i \label{E:E} \\[1ex]
H(P) & \le & \kappa q^{1+i} \label{E:H}
\end{eqnarray}
\end{lemma}

\noindent {\em Proof. \ }  In view of  (\ref{sv1}) and the fact that $j \le i$, it follows that
$$|E_P|\le|A_P|+(\kappa-1)|B_P|\le q^i+(\kappa-1)q^j\le \kappa q^i,$$
and hence
$$H(P)=\max\{q|E_P|,3\kappa\sqrt{cq|B_P|}\}\le\max\{\kappa q^{1+i},3\kappa c^{\frac{1}{2}}q^{\frac{1+j}{2}}\}=\kappa q^{1+i}.~\vspace*{-2ex}$$
 \hfill $ \boxtimes $

\vspace*{2ex}

In the above, to each $P\in\cP$  we have attached a rational line  $L_P $ with coefficients  $A_P,B_P,C_P$ and the function $F_P$ and the quantity $E_P$.  For ease of notation and the sake of  clarity,  we shall drop the subscript~$P$ if there is no ambiguity or confusion caused.

\vspace*{2ex}

\begin{lemma}\label{L:length}
If $P\in\cP^*$ (resp. $P\in\cP\setminus\cP^*$), then $\Theta(P)$ is contained in one (resp. at most two) open interval(s) of length at most $$\frac{42\kappa^3c}{H(P)}  \, . $$
\end{lemma}

\noindent{\em Proof. \ }   {\em Case (1).}  Suppose $P\in\cP^*$. Then $B\ne 0$. For any $x\in \Theta(P)$, it follows from Lemma \ref{L:Theta} and \eqref{E:cP*} that
\begin{eqnarray*}
\frac{3\kappa c}{q} & > & \Big|F(x)-F\Big(\frac{p}{q}\Big) \Big| \ = \ \Big|E\Big(x-\frac{p}{q}\Big)+\frac{F''(\xi)}{2}\Big(x-\frac{p}{q}\Big)^2\Big|  \\
& \ge & \frac{|B|}{2\kappa}\Big|x-\frac{p}{q}\Big|^2-3\kappa\sqrt{\frac{c|B|}{q}}\Big|x-\frac{p}{q}\Big| \, ,
\end{eqnarray*}
where $\xi\in I$. This implies that
$$\Big|x-\frac{p}{q}\Big|  \ <  \  7\kappa^2\sqrt{\frac{c}{q|B|}}  \ = \ \frac{21\kappa^3c}{H(P)} \, .$$
The upshot of this is that   $\Theta(P)$ is contained in an open interval of length $42\kappa^3c/H(P)$.

{\em Case (2).}  Suppose $P\in\cP\setminus\cP^*$. If $B=0$, then $E=A\ne0$ and $F(x)=Ex+C$. It follows that  $\Theta(P)$ is an open interval of length
$$\frac{4\kappa c}{q|E|} \ = \ \frac{4\kappa c}{H(P)} \ < \ \frac{42\kappa^3c}{H(P)}.$$

Now suppose $B\ne0$. Then, by \eqref{E:cP*}, we have that $E\ne0$. Consider the closed interval
$$\Omega:=\Big\{x\in I:\Big|x-\frac{p}{q}\Big|\le\frac{|E|}{\kappa|B|}\Big\}.$$
We first prove that $\Theta(P)\cap\Omega$ is contained in an open interval of length $42\kappa^3c/H(P)$. If $x\in\Theta(P)\cap\Omega$, then it follows from Lemma \ref{L:Theta} that
\begin{eqnarray*}
\frac{3\kappa c}{q} & > & \Big|F(x)-F\Big(\frac{p}{q}\Big)\Big|  \ = \ \Big|E\Big(x-\frac{p}{q}\Big)+\frac{1}{2}B f''(\xi)\Big(x-\frac{p}{q}\Big)^2\Big| \\[2ex]
& \ge &  \Big(|E|-\frac{\kappa|B|}{2} \Big|x-\frac{p}{q}\Big|\Big)\Big|x-\frac{p}{q}\Big| \ \ge  \ \frac{|E|}{2}\Big|x-\frac{p}{q}\Big|.
\end{eqnarray*}
This implies that
\begin{equation} \label{jin}
\Big|x-\frac{p}{q}\Big| \ < \ \frac{6\kappa c}{q|E|} \ = \ \frac{6\kappa c}{H(P)}.
\end{equation}
Hence $\Theta(P)\cap\Omega$ is contained in an open interval of length $12\kappa c/H(P)<42\kappa^3c/H(P)$. In particular, this implies the desired statement if $\Theta(P)\subset\Omega$.

Suppose that $\Theta(P)\not\subset\Omega$. Then the following three observations imply that $\Theta(P)\setminus\Omega$ is a connected interval.
\begin{itemize}
  \item $ |F(\frac{p}{q})|<\kappa c/q$ by Lemma \ref{L:Theta}.
  \item If $x_0$ is an end point of $\Omega$ and is contained in the interior of $I$, then $|F(x_0)| \ge 2\kappa c/q$. To see this, note that the assumption on $x_0$ and \eqref{E:cP*} imply that $$\Big|x_0-\frac{p}{q}\Big| \ = \ \frac{|E|}{\kappa|B|} \ > \ \frac{6\kappa c}{q|E|},$$ which together  with \eqref{jin} implies that $x_0\notin\Theta(P)$. Hence the desired inequality follows from \eqref{E:Theta}.
  \item The function $F:I\to\RR$ is either convex or concave.
\end{itemize}
We next claim that
\begin{equation} \label{sleep}
|F'(\xi)|\ge\frac{|E|}{2\kappa}   \qquad {\rm  for \ any  } \quad  \xi\in\Theta(P)\setminus\Omega .
\end{equation}
Assuming this claim for the moment, it follows that for any $x_1,x_2\in\Theta(P)\setminus\Omega$
$$\frac{4\kappa c}{q}  \  >  \ |F(x_1)-F(x_2)| \ = \ |F'(\xi)| \ |x_1-x_2|  \ \ge  \ \frac{|E|}{2\kappa}|x_1-x_2|,$$
where $\xi\in\Theta(P)\setminus\Omega$. Thus
$$|x_1-x_2| \ < \ \frac{8\kappa^2c}{q|E|} \ = \ \frac{8\kappa^2c}{H(P)}.$$
Hence $\Theta(P)\setminus\Omega$ is contained in an open interval of length $8\kappa^2c/H(P)<42\kappa^3c/H(P)$ and thereby completes  the proof of the lemma modulo (\ref{sleep}).

We now prove (\ref{sleep}) in the instance that $f''>0$ and $B>0$. The other cases are similar and left to the reader. Let $\xi\in\Theta(P)\setminus\Omega$, and consider the functions $G: \RR \to\RR$ and $H : \RR \to\RR$ given by
$$G(x): =\frac{\kappa B}{2}\Big(x-\frac{p}{q}\Big)^2+E \, \Big(x-\frac{p}{q}\Big)+F\Big(\frac{p}{q}\Big)$$
and
$$H(x):=\frac{B}{2\kappa}(x-\xi)^2+F'(\xi)(x-\xi)+F(\xi).$$
It follows that  $G(\frac{p}{q})=F(\frac{p}{q})$, $G'(\frac{p}{q})=F'(\frac{p}{q})$, $H(\xi)=F(\xi)$, and $H'(\xi)=F'(\xi)$. Moreover, if $x\in I$, then
\begin{eqnarray*}
G''(x)-F''(x) & = & B(\kappa-f''(x)) \ \ge \ 0, \\
H''(x)-F''(x) & = & B(\kappa^{-1}-f''(x)) \ \le \ 0.
\end{eqnarray*}
Thus, for $x\in I$ we have that $G(x)\ge F(x)\ge H(x)$. It follows that for any $x\in I$,
\begin{equation}\label{E:G-min}
G(x) \ \ge \ H(x) \ \ge \ -\frac{\kappa F'(\xi)^2}{2B}+F(\xi) \ \ge \ -\frac{\kappa F'(\xi)^2}{2B}-\frac{2\kappa c}{q}.
\end{equation}
Note that $F''=Bf''>0$. So the end point of $\Omega$ that is contained in the interval with end points $\frac{p}{q}$ and $\xi$ is equal to $\frac{p}{q}-\frac{E}{\kappa B}$. In particular, we have that $\frac{p}{q}-\frac{E}{\kappa B}\in I$. Thus \eqref{E:G-min} implies that
$$G\Big(\frac{p}{q}-\frac{E}{\kappa B}\Big) \ \ge \ -\frac{\kappa F'(\xi)^2}{2B}-\frac{2\kappa c}{q}.$$
On the other hand, we have that
$$G\Big(\frac{p}{q}-\frac{E}{\kappa B}\Big) \ = \ -\frac{E^2}{2\kappa B}+F\Big(\frac{p}{q}\Big) \ \le \ -\frac{E^2}{2\kappa B}+\frac{\kappa c}{q}.$$
On combining the previous two inequalities,  we find that
$$F'(\xi)^2 \ \ge \ \frac{E^2}{\kappa^2}-\frac{6cB}{q} \ \ge \ \frac{E^2}{\kappa^2}-\frac{2E^2}{3\kappa^2} \ = \ \frac{E^2}{3\kappa^2}$$
and  (\ref{sleep})  follows. This completes the proof of Lemma \ref{L:length}.
~\vspace*{-2ex}

 \hfill $ \boxtimes $

\section{Proof of Proposition \ref{P:main}  \label{ohohoh}}

For $n\ge1$, let
$$H_n:=42\kappa^3cl^{-1}R^n$$
and
\begin{equation}\label{E:P_n}\cP_n: =\Big\{P=\Big(\frac{p}{q},\frac{r}{q}\Big)\in\cP: H_n\le H(P)< H_{n+1}\Big\}.\end{equation}
Note that if $P\in\cP_n$, then, by \eqref{E:H}, we have that
\begin{equation}\label{missit}\kappa q^{1+i} \ge H_n.\end{equation}
Next let
\begin{equation}\label{E:cPn0}
\cP_{n,0}:=\cP_n\cap\cP^*
\end{equation}
and
\begin{equation}\label{E:cPnk}
\cP_{n,k}:=\{P\in\cP_n\setminus\cP^*:H_nR^{\lambda_{k-1}}\le \kappa q^{1+i}\le H_nR^{\lambda_k}\}   \quad {\rm for }  \quad 1\le k\le n\,,
\end{equation}
where $\lambda_k$ are defined by \eqref{vb201}.

\begin{lemma}\label{L:partition}  With $\cP_n$ and $\cP_{n,k}$ as above, we have that
$$\cP=\bigcup_{n=1}^\infty\cP_n  \quad  { and }  \quad  \cP_n=\bigcup_{k=0}^{n}\cP_{n,k} \, .  $$
\end{lemma}

\noindent{\em Proof. \ }   It is easily verified via \eqref{E:c} that $H(P)\ge3\kappa c^{\frac{1}{2}}$ for any $P\in\cP$, and that $H_1=42\kappa^3cl^{-1}R\le3\kappa c^{\frac{1}{2}}$.
Hence, $\cP=\bigcup_{n=1}^\infty\cP_n$.

Since for $P\in\cP_n\setminus\cP^*$ we have
$$q \ = \ \frac{H(P)^2}{qE^2} \ \le \ \frac{H_{n+1}^2}{9\kappa^2c} \ = \ \frac{H_n^2R^2}{9\kappa^2c},$$
it follows from \eqref{E:mu} that
\begin{eqnarray*}
\frac{\kappa q^{1+i}}{H_nR^{\lambda_n}} & \le & \kappa\frac{q^2}{H_n}R^{-\lambda_n} \ \le \ \kappa\frac{H_n^3R^4}{81\kappa^4c^2}R^{-\lambda_n}\\
& = & \frac{42^3}{81}\kappa^6cl^{-3}R^{-(\frac{2(1+i)}{j}-3)n+4-\mu}\\
& \le & (10\kappa^2l^{-1}R^{2-\frac{\mu}{3}})^3 \ \le \ 1.
\end{eqnarray*}
This together with \eqref{missit} implies that $\cP_n=\bigcup_{k=0}^{n}\cP_{n,k}$.
~\vspace*{-1ex}
 \hfill $ \boxtimes $

\vspace*{3ex}

We claim that the partition of $\cP$  given by Lemma \ref{L:partition} satisfies the requirement of Proposition \ref{P:main}. In other words, $\cP$ gives rise to a tree $\cS$  as described in \S\ref{winstrat} that contains  an $([R]-10)$-regular subtree. Recall that $\cS$ is itself a subtree of an $[R]$-regular rooted tree $\cT$.  The key towards establishing the claim is the following lemma and its corollary. For $k\ge0$, we let  $$k^+:=\max\{k,1\}  \, .  $$
The following lemma contains a crucial property of the lines $L_P$ defined by \eqref{vb210}.

\begin{lemma}\label{L:const}
For any $n\ge1$, $0\le k\le n$ and $\tau\in\cS_{n-k^+}$, the map $P\mt L_P$ is constant on $$\cP_{n,k}(\tau):=\{P\in\cP_{n,k}:\cI(\tau)\cap\Delta(P)\ne\emptyset\}.$$
\end{lemma}

We postpone the proof  for the moment and continue by stating an important consequence of the lemma.

\begin{corollary}\label{C:est}
For any $n\ge1$, $0\le k\le n$ and $\tau\in\cS_{n-k^+}$, we have
$$\#\Big\{\tau'\in\cT_n:\cI(\tau')  \ \cap\bigcup_{P\in\cP_{n,k}(\tau)}\Delta(P)\ne\emptyset\Big\} \ \le \
\begin{cases}
2, &\text{if $k=0$,}\\
4, &\text{if $k\ge1$.}
\end{cases}$$
\end{corollary}

\noindent{\em Proof. \ }
We may assume that $\cP_{n,k}(\tau)\ne\emptyset$. Let $P_0=(\frac{p_0}{q_0},\frac{r_0}{q_0})\in\cP_{n,k}(\tau)$ be such that $q_0\le q$ for any $(\frac{p}{q},\frac{r}{q})\in\cP_{n,k}(\tau)$. By Lemma~\ref{L:const}, for any $P\in\cP_{n,k}(\tau)$ we have $L_P=L_{P_0}$ and so $\Theta(P)\subset\Theta(P_0)$. Thus it follows from Lemma~\ref{L:Theta} that
$$\bigcup_{P\in\cP_{n,k}(\tau)}\Delta(P) \  \subset \bigcup_{P\in\cP_{n,k}(\tau)}\Theta(P) \ \subset \ \Theta(P_0).$$
Hence
$$\Big\{\tau'\in\cT_n:\cI(\tau') \ \cap\bigcup_{P\in\cP_{n,k}(\tau)}\Delta(P)\ne\emptyset\Big\}
\ \subset\ \{\tau'\in\cT_n:\cI(\tau') \ \cap   \  \Theta(P_0)\ne\emptyset\}.$$
By Lemma \ref{L:length}, if $k=0$ (resp. $k\ge1$), then $\Theta(P_0)$ is contained in one (resp. at most two) open interval(s) of length at most
$$\frac{42\kappa^3c}{H(P_0)} \ \stackrel{\eqref{E:P_n}}{\le} \ \frac{42\kappa^3c}{H_n} \ = \ lR^{-n}.$$
Since the intervals $\{\cI(\tau'):\tau'\in\cT_n\}$ are of length $lR^{-n}$ and have mutually disjoint interiors, there can be at most $2$ (resp. $4$) of them that intersect $\Theta(P_0)$. This proves the corollary.
\vspace*{-2ex}

\hfill $ \boxtimes $

\vspace*{2ex}

We are now in the position to prove Proposition \ref{P:main}.   In view of Proposition \ref{propJA},  it suffices to prove that the intersection of $\cS$ with every $11$-regular subtree of $\cT$ is infinite.     Let $\cR\subset\cT$ be an $11$-regular subtree and let
$$\cR':=\cR\cap\cS   \qquad {\rm   and }  \qquad a_n:=\#\cR'_n  \quad (n \ge 0) \, , $$
where $\cR'_n$ is the $n$'th level of the tree $\cR'$.
Then $a_0=1$. We prove that  $\cR'$ is  infinite by showing that
\begin{equation}\label{E:infinity}
a_n>3a_{n-1}    \quad (n \ge 1) \, .
\end{equation}
We use induction. For $ n \ge 1$, let
$$\cU_n:=\Big\{\tau\in\cT_{\suc}(\cR'_{n-1}):\cI(\tau) \ \cap\bigcup_{P\in\cP_n}\Delta(P)\ne\emptyset\Big\}.$$
Then
$$
\cR'_n
 \ = \ \Big\{\tau\in\cR_{\suc}(\cR'_{n-1}):\cI(\tau) \ \cap\bigcup_{P\in\cP_n}\Delta(P)=\emptyset\Big\}
 \ = \ \cR_{\suc}(\cR'_{n-1})\setminus\cU_n.
$$
It follows that
\begin{equation}\label{E:infinity1}
a_n\ge 11a_{n-1}-\#\cU_n.
\end{equation}
On the other hand,
\begin{eqnarray*}
\cU_n & = & \bigcup_{k=0}^n\Big\{\tau\in\cT_{\suc}(\cR'_{n-1}):\cI(\tau) \ \cap\bigcup_{P\in\cP_{n,k}}\Delta(P)\ne\emptyset\Big\}\\
& \subset & \bigcup_{k=0}^n\Big\{\tau\in\cT_n:\tau\prec\cR'_{n-k^+},\cI(\tau) \ \cap\bigcup_{P\in\cP_{n,k}}\Delta(P)\ne\emptyset\Big\}\\
& = & \bigcup_{k=0}^n\bigcup_{~ \ \tau'\in\cR'_{n-k^+}}\!\!\! \Big\{\tau\in\cT_n:\tau\prec\tau',\cI(\tau) \ \cap\bigcup_{P\in\cP_{n,k}}\Delta(P)\ne\emptyset\Big\}\\
& \subset & \bigcup_{k=0}^n\bigcup_{~ \ \tau'\in\cR'_{n-k^+}} \!\!\! \Big\{\tau\in\cT_n:\cI(\tau) \ \cap\bigcup_{P\in\cP_{n,k}(\tau')}\Delta(P)\ne\emptyset\Big\}.
\end{eqnarray*}
Thus,  Corollary \ref{C:est} implies that
\begin{equation}\label{E:infinity2}
\#\cU_n\le 2a_{n-1}+\sum_{k=1}^n4a_{n-k}.
\end{equation}
On combining  \eqref{E:infinity1} and \eqref{E:infinity2}, we obtain that
\begin{equation}\label{E:infinity3}
a_n\ge 9a_{n-1}-\sum_{k=1}^n4a_{n-k}.
\end{equation}
With $n=1$ in \eqref{E:infinity3}, we  find that $a_1\ge 5$. Hence,  \eqref{E:infinity} holds for $n=1$. Now  assume $n\ge2$ and  that \eqref{E:infinity} holds with $n$  replaced by $1,\ldots,n-1$.
Then for any $1\le k\le n$, we have that  $$a_{n-k}\le 3^{-k+1}a_{n-1}.$$ Substituting this into \eqref{E:infinity3}, gives that
$$a_n\ge9a_{n-1}-4a_{n-1}\sum_{k=1}^n3^{-k+1}>3a_{n-1}.$$
This completes the induction step and thus establishes  \eqref{E:infinity}. In turn this completes the proof of Proposition \ref{P:main} modulo the truth of Lemma \ref{L:const}.
\vspace*{-2ex}

\hfill $ \boxtimes $

\vspace*{2ex}

\subsection{Proof of Lemma \ref{L:const}}

To begin with we prove the following result.

\begin{lemma}
Let $n,k $ and $ \tau$ be as in Lemma \ref{L:const}, and let $P_1:=(\frac{p_1}{q_1},\frac{r_1}{q_1}),P_2:=(\frac{p_2}{q_2},\frac{r_2}{q_2})\in\cP_{n,k}(\tau)$. Denote $A_s=A_{P_s}$, $B_s=B_{P_s}$, $C_s=C_{P_s}$ and $F_s=F_{P_s}$, $s=1,2$.
\begin{itemize}
\item[(1)]  If $k=0$ and $q_1\le q_2$, then
\begin{equation}
A_2p_1+B_2r_1+C_2q_1 \ = \ 0 \label{this1} \,
\end{equation}
\begin{equation} \label{this2}
|A_1B_2-A_2B_1| \ < \ q_1  \, .
\end{equation}
\item[(2)] If $k=1$, then
\begin{equation}\label{this3}
A_2p_1+B_2r_1+C_2q_1 \ = \ 0 \ = \  A_1p_2+B_1r_2+C_1q_2 \, .
\end{equation}
\item[(3)] If $k\ge2$, then
\begin{eqnarray}
\Big|A_1F_2\Big(\frac{p_1}{q_1}\Big)-A_2F_1\Big(\frac{p_1}{q_1}\Big)\Big| & \le & \frac{c}{2}H_n^{-\frac{j}{1+i}}, \label{this4} \\[1ex]
\Big|B_1F_2\Big(\frac{p_1}{q_1}\Big)-B_2F_1\Big(\frac{p_1}{q_1}\Big)\Big| & \le & \frac{c}{2}H_n^{-\frac{i}{1+i}}, \label{this5} \\[1ex]
|A_1B_2-A_2B_1| & \le & \frac{1}{3\kappa}H_n^{\frac{1}{1+i}}R^{-k}.  \label{this6}
\end{eqnarray}
\end{itemize}
\end{lemma}

\noindent{\em Proof. \ }
We first prove that for any $0\le k\le n$,
\begin{eqnarray}
\Big|\frac{p_1}{q_1}-\frac{p_2}{q_2}\Big| & \le & 2lR^{-n+k^+},\label{E:pq} \\[1ex]
\max \Big\{ q_2\Big|F_2\Big(\frac{p_1}{q_1}\Big)\Big|, \, q_1\Big|F_1\Big(\frac{p_2}{q_2}\Big)\Big|  \Big\} & \le & 477\kappa^5cR^{2k^++2}, \label{sva}\\[1ex]
|A_2p_1+B_2r_1+C_2q_1| & \le & \frac{1}{2}R^{-4-\lambda_1}\Big(\frac{q_1}{q_2}R^{2k^++2}+\frac{|B_2|}{q_1^j}\Big),  \label{svb}\\[1ex]
|A_1p_2+B_1r_2+C_1q_2| & \le & \frac{1}{2}R^{-4-\lambda_1}\Big(\frac{q_2}{q_1}R^{2k^++2}+\frac{|B_1|}{q_2^j}\Big),  \label{svba}\\[1ex]
|A_1B_2-A_2B_1| & \le & |E_1B_2|+|E_2B_1|+2|B_1B_2|\kappa lR^{-n+k^+}\,, \label{svc}
\end{eqnarray}
where $E_1=E_{P_1}$ and $E_2=E_{P_2}$.

To establish \eqref{E:pq}, note that for $s=1,2$ there exists $x_s\in\cI(\tau)\cap\Delta(P_s)\ne\emptyset$ and that $|x_1-x_2|\le \rho(\cI(\tau))=lR^{-n+k^+}$. Then, it follows from \eqref{E:Delta} that
\begin{eqnarray*} \Big|\frac{p_1}{q_1}-\frac{p_2}{q_2}\Big| & \le & \Big|\frac{p_1}{q_1}-x_1\Big|+
|x_1-x_2|+\Big|x_2-\frac{p_2}{q_2}\Big|  \\[2ex]
&\le &  \frac{c}{q_1^{1+i}}+lR^{-n+k^+}+\frac{c}{q_2^{1+i}}
\ \stackrel{\eqref{missit}}{\le} \ \frac{2\kappa c}{H_n}+lR^{-n+k^+}
\\[1ex] & = &   \ \frac{1}{21\kappa^2}lR^{-n}+lR^{-n+k^+}
\ \le   \  2lR^{-n+k^+}.\end{eqnarray*}
Regarding \eqref{sva}, we expand $F_2(\frac{p_1}{q_1})$ using Taylor's formula at the point $\frac{p_2}{q_2}$ and estimate as follows
\begin{eqnarray*}
q_2\Big|F_2\Big(\frac{p_1}{q_1}\Big)\Big| & \le & q_2\Big|F_2\Big(\frac{p_2}{q_2}\Big)\Big|+q_2 |E_2| \ \Big|\frac{p_1}{q_1}-\frac{p_2}{q_2}\Big|+\frac{\kappa q_2|B_2|}{2} \ \Big|\frac{p_1}{q_1}-\frac{p_2}{q_2}\Big|^2\\[1ex]
& \le &  \kappa c+H(P_2)\Big|\frac{p_1}{q_1}-\frac{p_2}{q_2}\Big|+\frac{H(P_2)^2}{18\kappa c}\Big|\frac{p_1}{q_1}-\frac{p_2}{q_2}\Big|^2\\[1ex]
& \stackrel{\eqref{E:pq}}{\le}  & \kappa c+H_{n+1}\cdot2lR^{-n+k^+}+\frac{H_{n+1}^2}{18\kappa c}\cdot4l^2R^{-2n+2k^+}\\[1ex]
& =  & \kappa c+84\kappa^3cR^{k^++1}+392\kappa^5cR^{2k^++2}\\[1ex]
& \le  & 477\kappa^5cR^{2k^++2}.
\end{eqnarray*}
The above argument can be trivially  modified  to show that the same upper bound  is valid for   $q_1|F_1(\frac{p_2}{q_2})|$.  Turning our attention to  \eqref{svb}, using \eqref{E:cP} we estimate as follows
\begin{eqnarray*}
|A_2p_1+B_2r_1+C_2q_1|
& = & q_1\Big|F_2\Big(\frac{p_1}{q_1}\Big)-B_2\Big(f\Big(\frac{p_1}{q_1}\Big)-\frac{r_1}{q_1}\Big)\Big|\\[1ex]
&\le&q_1\Big|F_2\Big(\frac{p_1}{q_1}\Big)\Big|+q_1|B_2|\Big|f\Big(\frac{p_1}{q_1}\Big)-\frac{r_1}{q_1}\Big|\\[1ex]
&\stackrel{\eqref{sva}}{\le}& 477\kappa^5c\frac{q_1}{q_2}R^{2k^++2}+|B_2|\frac{\kappa c}{q_1^j}\\[1ex]
&\le& 477\kappa^5c\Big(\frac{q_1}{q_2}R^{2k^++2}+\frac{|B_2|}{q_1^j}\Big)\\[1ex]
&\stackrel{\eqref{E:c}}{\le}& \frac{1}{2}R^{-4-\lambda_1}\Big(\frac{q_1}{q_2}R^{2k^++2}+\frac{|B_2|}{q_1^j}\Big).
\end{eqnarray*}
The same argument with obvious modifications yields \eqref{svba}.
Finally, regarding \eqref{svc}, using the definition \eqref{vb2003} for $E_1$ and $E_2$, we have that
\begin{eqnarray*}
|A_1B_2-A_2B_1| & = & \Big|\Big(E_1-f'\Big(\frac{p_1}{q_1}\Big)B_1\Big)B_2-\Big(E_2-f'\Big(\frac{p_2}{q_2}\Big)B_2\Big)B_1\Big|\\
& = & \Big|E_1B_2-E_2B_1+B_1B_2\Big(f'\Big(\frac{p_2}{q_2}\Big)-f'\Big(\frac{p_1}{q_1}\Big)\Big)\Big|\\
& \le  & |E_1B_2|+|E_2B_1|+|B_1B_2|\kappa\Big|\frac{p_1}{q_1}-\frac{p_2}{q_2}\Big|\\
& \stackrel{\eqref{E:pq}}{\le} &  |E_1B_2|+|E_2B_1|+2|B_1B_2|\kappa lR^{-n+k^+}.
\end{eqnarray*}
Having established  \eqref{E:pq}--\eqref{svc}, we are now in the position to prove the lemma.

\vspace*{2ex}

{\em Part (1). \ }  Suppose $k=0$ and $q_1\le q_2$. It follows from the definition \eqref{E:cPn0} for $\cP_{n,0}$ that
\begin{eqnarray*}
|A_2p_1+B_2r_1+C_2q_1| & \stackrel{\eqref{svb}}{\le}  & \frac{1}{2}R^{-4-\lambda_1} \Big(\frac{q_1}{q_2}R^4+\frac{|B_2|}{q_1^j}\Big)\\[1ex]
& \le  & \frac{1}{2}R^{-4-\lambda_1}\Big(R^4+\frac{q_2|B_2|}{q_1|B_1|}\Big) \\[1ex]
& =  & \frac{1}{2}R^{-4-\lambda_1}\Big(R^4+\frac{H(P_2)^2}{H(P_1)^2}\Big)\\[1ex]
& \le  & \frac{1}{2}R^{-4-\lambda_1}(R^4+R^2)\\[1ex]
& \le  &  R^{-\lambda_1} \ < \ 1.
\end{eqnarray*}
Since the left hand side of the above inequality is an integer, it follows that  $A_2p_1+B_2r_1+C_2q_1=0$.  Regarding \eqref{this2},  note that for $s=1,2$
$$|E_s| \ \le \ \frac{H(P_s)^2}{q_sH(P_s)} \ \le \ \frac{9\kappa^2c|B_s|}{H_n} \ \le \ \frac{lR^{-n}}{2\kappa}|B_s|$$
and that
$$\frac{|B_2|}{|B_1|} \ \le \ \frac{q_2|B_2|}{q_1|B_1|} \ = \ \frac{H(P_2)^2}{H(P_1)^2} \ \le \ \frac{H_{n+1}^2}{H_n^2} \ = \ R^2.$$
These inequalities together with  \eqref{E:l1} imply that
\begin{eqnarray*}
|A_1B_2-A_2B_1| & \stackrel{\eqref{svc}}{\le} & |E_1B_2|+|E_2B_1|+2|B_1B_2|\kappa lR^{-n+1}\\
& \le & |B_1B_2|\Big(\frac{lR^{-n}}{\kappa}+2\kappa lR^{-n+1}\Big)\\
& \le & B_1^2R^2\cdot3\kappa l \\
& \le & 3\kappa lR^2q_1^{2j} \ < \ q_1.
\end{eqnarray*}

{\em Part (2). \ }  Suppose $k=1$. It follows from \eqref{E:cPnk} that
$$\max\Big\{\frac{q_1}{q_2},\frac{q_2}{q_1}\Big\}\le R^{\frac{\lambda_1}{1+i}}<R^{\lambda_1}.$$ This implies that
\begin{eqnarray*}
|A_2p_1+B_2r_1+C_2q_1| & \stackrel{\eqref{svb}}{\le} & \frac{1}{2}R^{-4-\lambda_1}\Big(\frac{q_1}{q_2}R^4+\frac{|B_2|}{q_1^j}\Big) \\
&\le& \frac{1}{2}R^{-\lambda_1}\Big(\frac{q_1}{q_2}+\frac{q_2^j}{q_1^j}\Big)  \ < \ 1 \, .
\end{eqnarray*}
The left hand side is an integer and so  must be zero.   The same argument involving \eqref{svba} rather than \eqref{svb} shows that $|A_1p_2+B_1r_2+C_1q_2|=0$.

{\em Part (3). \ } Suppose $k\ge2$. We first prove that
\begin{eqnarray}
\max\Big\{\frac{q_1^i}{q_2},\frac{q_2^i}{q_1}\Big\} & \le & \frac{1}{10^3\kappa^5}H_n^{-\frac{j}{1+i}}R^{-2k-2} \label{i1} \\[1ex]
\max\Big\{\frac{q_1^j}{q_2},\frac{q_2^j}{q_1}\Big\} & \le & \frac{1}{10^3\kappa^5}H_n^{-\frac{i}{1+i}}R^{-2k-2}  \label{i2} \, .
\end{eqnarray}
It follows from \eqref{E:cPnk} that
$$\max\{q_1,q_2\} \ \le \ (\kappa^{-1}H_n)^{\frac{1}{1+i}}R^{\frac{\lambda_k}{1+i}}$$
and
$$\min\{q_1,q_2\} \ \ge \ (\kappa^{-1}H_n)^{\frac{1}{1+i}}R^{\frac{\lambda_{k-1}}{1+i}}.$$
In view of the fact that  $$\frac{j\lambda_k-\lambda_{k-1}}{1+i} \ \le \ \frac{i\lambda_k-\lambda_{k-1}}{1+i} \ = \ -\frac{j}{1+i}\mu+\frac{2}{j}-2k,$$
it follows from \eqref{E:mu} that
\begin{eqnarray*}
\max\Big\{\frac{q_1^i}{q_2},\frac{q_2^i}{q_1}\Big\} & \le & (\kappa^{-1}H_n)^{-\frac{j}{1+i}}R^{\frac{i\lambda_k-\lambda_{k-1}}{1+i}} \\
& \le & \kappa H_n^{-\frac{j}{1+i}}R^{-\frac{j}{1+i}\mu+\frac{2}{j}-2k}\\[1ex]
& \le & \frac{1}{10^3\kappa^5}H_n^{-\frac{j}{1+i}}R^{-2k-2}
\end{eqnarray*}
and that
\begin{eqnarray*}
\max\Big\{\frac{q_1^j}{q_2},\frac{q_2^j}{q_1}\Big\}& \le & (\kappa^{-1}H_n)^{-\frac{i}{1+i}}R^{\frac{j\lambda_k-\lambda_{k-1}}{1+i}}\\
& \le & \kappa H_n^{-\frac{i}{1+i}}R^{-\frac{j}{1+i}\mu+\frac{2}{j}-2k}\\[1ex]
& \le & \frac{1}{10^3\kappa^5}H_n^{-\frac{i}{1+i}}R^{-2k-2}.
\end{eqnarray*}
This establishes \eqref{i1} and \eqref{i2}. It now follows that
\begin{eqnarray*}
\Big|A_1F_2\Big(\frac{p_1}{q_1}\Big)-A_2F_1\Big(\frac{p_1}{q_1}\Big)\Big| & \le & |A_1|\Big|F_2\Big(\frac{p_1}{q_1}\Big)\Big|+|A_2|\Big|F_1\Big(\frac{p_1}{q_1}\Big)\Big|\\
& \stackrel{\eqref{sva}}{\le} & q_1^i\cdot\frac{1}{q_2}477\kappa^5cR^{2k+2}+q_2^i\frac{\kappa c}{q_1}\\
& \stackrel{\eqref{i1}}{\le} & \frac{1}{10^3\kappa^5}H_n^{-\frac{j}{1+i}}R^{-2k-2}(477\kappa^5cR^{2k+2}+\kappa c)\\
& \le & \frac{c}{2}H_n^{-\frac{j}{1+i}}
\end{eqnarray*}
and
\begin{eqnarray*}
\Big|B_1F_2\Big(\frac{p_1}{q_1}\Big)-B_2F_1\Big(\frac{p_1}{q_1}\Big)\Big| & \le & |B_1|\Big|F_2\Big(\frac{p_1}{q_1}\Big)\Big|+|B_2|\Big|F_1\Big(\frac{p_1}{q_1}\Big)\Big|\\
& \stackrel{\eqref{sva}}{\le} & q_1^j\cdot\frac{1}{q_2}477\kappa^5cR^{2k+2}+q_2^j\frac{\kappa c}{q_1}\\
& \stackrel{\eqref{i2}}{\le} & \frac{1}{10^3\kappa^5}H_n^{-\frac{i}{1+i}}R^{-2k-2}(477\kappa^5cR^{2k+2}+\kappa c) \\
& \le &\frac{c}{2}H_n^{-\frac{i}{1+i}}\,.
\end{eqnarray*}
Finally, since $P_s\not\in\cP^*$, we have that $q_sE_s^2\ge 9\kappa^2c|B_s|$ and $H(P_s)=q_s|E_s|$ for $s=1,2$. Then
$$|B_2| \ \le \ \frac{H(P_2)^2}{9\kappa^2cq_2}$$
and using \eqref{vb204} and the fact that $H_{n+1}=RH_n$ , we get that
\begin{eqnarray*}
|A_1B_2-A_2B_1| & \stackrel{\eqref{svc}}{\le} & |E_1B_2|+|E_2B_1|+2|B_1B_2|\kappa lR^{-n+k}\\[1ex]
& \le & \frac{H(P_1)}{q_1}q_2^j+\frac{H(P_2)}{q_2}q_1^j+2q_1^j\frac{H(P_2)^2}{9\kappa^2cq_2}\kappa lR^{-n+k}\\[1ex]
& \stackrel{\eqref{i2}}{\le} & \frac{1}{10^3\kappa^5}H_n^{-\frac{i}{1+i}}R^{-2k-2}\Big(2H_{n+1}+\frac{2H_{n+1}^2}{9\kappa^2c}\kappa lR^{-n+k}\Big)\\[1ex]
& \le & \frac{1}{500\kappa^5}H_n^{\frac{1}{1+i}}R^{-2k-2}(R+5\kappa^2R^{k+2})\\[1ex]
& \le & \frac{1}{3\kappa}H_n^{\frac{1}{1+i}}R^{-k}.
\end{eqnarray*}
This thereby completes the proof of Part (3) and thus the lemma.
\vspace*{-2ex}

\hfill $ \boxtimes $

\vspace*{3ex}

We now proceed  with the proof of Lemma \ref{L:const}.  Let $P_1=(\frac{p_1}{q_1},\frac{r_1}{q_1})$ and $P_2=(\frac{p_2}{q_2},\frac{r_2}{q_2})$ be distinct points in $\cP_{n,k}(\tau)$. We need to prove that $L_{P_1}=L_{P_2}$.  We consider three separate cases.
\vspace*{2ex}

{\em Case (1). \ } Suppose $k=0$. Without loss of generality, we assume that $q_1\le q_2$.
Then in view of \eqref{this1} we have that  $A_2p_1+B_2r_1+C_2q_1=0$. Hence $L_{P_2}$ passes through $P_1$.
We prove that $L_{P_1}=L_{P_2}$ by contradiction. Thus, assume that  $L_{P_1}\ne L_{P_2}$. Since  the two lines intersect at $P_1$, it follows that
$$
\frac{p_1}{q_1}   \  =  \  \frac{B_1C_2 - B_2C_1}{A_1B_2-A_2B_1}        \qquad {\rm  and } \qquad \frac{r_1}{q_1} \  = \  \frac{A_2C_1 - A_1C_2}{A_1B_2-A_2B_1}  \, .
$$
In particular, since $p_1,q_1,r_1 $ are co-prime, the non-zero integer $ A_1B_2-A_2B_1 $  is divisible by $q_1$.  Thus
$$q_1 \ \le \ |A_1B_2-A_2B_1| \ \stackrel{\eqref{this2}}{<} \  q_1 \, , $$
which is of course impossible.

\vspace*{2ex}

{\em Case (2). \ } Suppose $k=1$. Then in view of \eqref{this3} we have that  $A_2p_1+B_2r_1+C_2q_1=0= A_1p_2+B_1r_2+C_1q_2$. Hence,  $L_{P_2}$ passes through $P_1$ and  $L_{P_1}$ passes through $P_2$. By definition, $L_{P_2}$ passes through $P_2$ and  $L_{P_1}$ passes through $P_1$.  The upshot is that both lines pass through the points $P_1$ and $P_2$,  and so we must have that  $L_{P_1}=L_{P_2}$.

\vspace*{2ex}

{\em Case (3). \ } Suppose $k\ge2$. We prove that $L_{P_1}=L_{P_2}$ by contradiction. Thus, assume that $L_{P_1}\ne L_{P_2}$. We first consider the case where $L_{P_1}$ is parallel to $L_{P_2}$. Then, it is easily verifed that
$$A_1F_2\Big(\frac{p_1}{q_1}\Big)-A_2F_1\Big(\frac{p_1}{q_1}\Big) \ = \ A_1C_2-A_2C_1$$
and
$$B_1F_2\Big(\frac{p_1}{q_1}\Big)-B_2F_1\Big(\frac{p_1}{q_1}\Big) \ = \ B_1C_2-B_2C_1 \, . $$
Since $H_n\ge c$, it follows via \eqref{this4} and  \eqref{this5} that
\begin{eqnarray*}
1 & \le & |A_1C_2-A_2C_1|+|B_1C_2-B_2C_1|\\[1ex]
 & \le & \Big|A_1F_2\Big(\frac{p_1}{q_1}\Big)-A_2F_1\Big(\frac{p_1}{q_1}\Big)\Big|+\Big|B_1F_2\Big(\frac{p_1}{q_1}\Big)-B_2F_1\Big(\frac{p_1}{q_1}\Big)\Big|\\[1ex]
 & \le & \frac{c}{2}H_n^{-\frac{j}{1+i}}+\frac{c}{2}H_n^{-\frac{i}{1+i}}\\[1ex]
 & \le & \frac{c}{2}(c^{-\frac{j}{1+i}}+c^{-\frac{i}{1+i}})\\[1ex]
 & = & \frac{1}{2}(c^{\frac{2i}{1+i}}+c^{\frac{1}{1+i}}) \ < \ 1
\end{eqnarray*}
which is of course impossible.

Now suppose $L_{P_1}$ is not parallel to $L_{P_2}$. Let $P_0=(\frac{p_0}{q_0},\frac{r_0}{q_0})\in\QQ^2$ be the point of  intersection of $L_{P_1}$ and $L_{P_2}$. Then it follows that the non-zero integer $A_1B_2-A_2B_1$ is  divisible by $q_0$ and so
\begin{equation} \label{1sv}  q_0 \ \le \ |A_1B_2-A_2B_1|.\end{equation}
We first prove that  $\Delta(P_1)\subset\Delta(P_0)$ and that $P_0\in\cP$.
It is easily verified that
$$\begin{pmatrix}
A_1 & B_1 \\ A_2 & B_2
\end{pmatrix}
\begin{pmatrix}
\frac{p_1}{q_1}-\frac{p_0}{q_0} \\ f(\frac{p_1}{q_1})-\frac{r_0}{q_0}
\end{pmatrix}
 \ = \ \begin{pmatrix}
F_1(\frac{p_1}{q_1}) \\ F_2(\frac{p_1}{q_1})
\end{pmatrix}.$$
Hence, on using Cramer's rule, we obtain that
\begin{equation}  \label{2sv} |A_1B_2-A_2B_1|\Big|\frac{p_1}{q_1}-\frac{p_0}{q_0}\Big| \ = \
\Big|B_1F_2\Big(\frac{p_1}{q_1}\Big)-B_2F_1\Big(\frac{p_1}{q_1}\Big)\Big|
 \ \stackrel{\eqref{this5}}{\le} \  \frac{c}{2}H_n^{-\frac{i}{1+i}},\end{equation}
and
\begin{equation}  \label{3sv} |A_1B_2-A_2B_1|\Big|f\Big(\frac{p_1}{q_1}\Big)-\frac{r_0}{q_0}\Big| \ = \
\Big|A_1F_2\Big(\frac{p_1}{q_1}\Big)-A_2F_1\Big(\frac{p_1}{q_1}\Big)\Big|
 \ \stackrel{\eqref{this4}}{\le} \ \frac{c}{2}H_n^{-\frac{j}{1+i}}.\end{equation}
If $x\in\Delta(P_1)$, then \eqref{this6} and  \eqref{2sv} imply that
\begin{eqnarray*}
q_0^{1+i}\Big|x-\frac{p_0}{q_0}\Big| & \stackrel{\eqref{1sv}}{\le} & |A_1B_2-A_2B_1|^{1+i}\Big|x-\frac{p_1}{q_1}\Big|+|A_1B_2-A_2B_1|^{1+i}\Big|\frac{p_1}{q_1}-\frac{p_0}{q_0}\Big|\\
& \le &  \frac{H_n}{3\kappa}\cdot\frac{c}{q_1^{1+i}}+H_n^{\frac{i}{1+i}}\cdot \frac{c}{2}H_n^{-\frac{i}{1+i}}\\
& \stackrel{\eqref{missit}}{\le} & \frac{c}{3}+\frac{c}{2} \ < \ c.
\end{eqnarray*}
Thus $x\in\Delta(P_0)$  and the upshot is that $\Delta(P_1)\subset\Delta(P_0)$. In particular,
\begin{equation} \label{this}
\Big|\frac{p_1}{q_1}-\frac{p_0}{q_0}\Big| \ < \ \frac{c}{q_0^{1+i}}.\end{equation}
Since $\cA_0\cap\Delta(P_1)\ne\emptyset$, there exists $x\in\cA_0$ such that $|x-\frac{p_1}{q_1}|<c/q_1^{1+i}$, and hence it follows that
$$\Big|x-\frac{p_0}{q_0}\Big| \ \le \ \Big|x-\frac{p_1}{q_1}\Big|+ \Big| \frac{p_1}{q_1}-\frac{p_0}{q_0}\Big| \ \le \ 2c \ \le \ \frac{l}{2}.$$
This implies that $\frac{p_0}{q_0}\in\cB_0\subset I$.  Also  note that by  \eqref{this6} and  \eqref{3sv},  we have that
\begin{eqnarray*}
q_0^{1+j}\Big|f\Big(\frac{p_0}{q_0}\Big)-\frac{r_0}{q_0}\Big|
& \stackrel{\eqref{1sv}}{\le} & q_0^{1+j}\Big|f\Big(\frac{p_0}{q_0}\Big)-f\Big(\frac{p_1}{q_1}\Big)\Big|+|A_1B_2-A_2B_1|^{1+j}\Big|f\Big(\frac{p_1}{q_1}\Big)-\frac{r_0}{q_0}\Big|\\[2ex]
& \le & q_0^{1+i}(\kappa-1)\Big|\frac{p_0}{q_0}-\frac{p_1}{q_1}\Big|+H_n^{\frac{j}{1+i}}\cdot\frac{c}{2}H_n^{-\frac{j}{1+i}}\\[2ex]
& \stackrel{\eqref{this}}{\le} & (\kappa-1)c+\frac{c}{2}\\[2ex]
& < & \kappa c.
\end{eqnarray*}
Thus $P_0\in\cP$ and so there exists a  unique integer $n_0\ge1$  such that $P_0\in\cP_{n_0}$.   Suppose for the moment that  $n_0\le n-k$.  Then there exists $\tau'\in\cS_{n_0}$ such that $\tau\prec\tau'$, and hence $$\cI(\tau)\cap\Delta(P_1) \ \subset  \  \cI(\tau')\cap\Delta(P_0)=\emptyset.$$ This contradicts the fact that $P_1\in\cP_{n,k}(\tau)$.  Thus
$$n_0\ge n-k+1 \, , $$
and so
$$H(P_0) \ \ge \  H_{n_0} \ \ge \  H_{n-k+1}.$$
On the other hand, we have that
$$H(P_0) \ \le \  \kappa q_0^{1+i} \ \stackrel{\eqref{1sv}}{\le} \  \kappa|A_1B_2-A_2B_1|^{1+i} \ \stackrel{\eqref{this6}}{\le} \  H_nR^{-k} \ = \ H_{n-k}  \, . $$
This contradicts the above lower bound for $H(P_0)$ and so completes the proof of Case (3) and indeed the lemma.
\vspace*{-2ex}

\hfill $ \boxtimes $

\vspace*{2ex}

\section{The inhomogeneous case: establishing  Theorem  \ref{thm1} \label{inhThm1.2} }

Theorem \ref{thm1} is easily deduced from the following  statement.

\begin{theorem}\label{T:inhomogeneous}
Let  $(i,j)$ be a pair of real numbers
satisfying   $0<j\le i<1$ and $i+j=1$.   Let  $I\subset\RR$ be a compact interval and  $f \in C^{(2)}(I) $ such  that  $ f''(x) \neq  0 $ for all $ x \in I $.   Then, for any $\bm\theta=(\gamma,\delta)\in\RR^2$ we have that
$\Bad^f_{\bm\theta}(i,j)$
is a $1/2$-winning subset of $I$.
\end{theorem}

\noindent We have already established  the homogeneous case ($\gamma=\delta=0$) of the Theorem~\ref{T:inhomogeneous};  namely Theorem \ref{T:main}.  With reference to \S\ref{winstrat},  the crux of the `homogeneous' proof involved constructing a partition $\cP_{n } $   ($ n \ge 1$)  of $\cP$  (given by  Lemma \ref{L:partition})  such that the subtree $\cS$ of an $[R]$-regular rooted tree $\cT$ has an $([R]-10)$-regular subtree $\cS'$ -- the substance of Proposition \ref{P:main}.  Recall,  by construction we have that
$$\cI(\tau) \ \cap\bigcup_{P\in\cP_n}\Delta(P)=\emptyset, \qquad \forall  \   n\ge1  \quad { \rm and } \quad  \tau\in\cS'_n.$$
To prove  Theorem~\ref{T:inhomogeneous},  the idea is to merge the inhomogeneous constraints into the homogeneous construction.  More precisely, we show that $\cS'$ has an $([R]-12)$-regular subtree $\cQ'$ that incorporates the inhomogeneous constraints.   With this in mind, let
$$c':=\frac{1}{10}cR^{-2}$$
where $c$ is defined in \eqref{E:c}, and let
$$\sV:=\Big\{(p,r,q)\in\ZZ^2\times\NN : \frac{p+\gamma}{q}\in I, \ \Big|f\Big(\frac{p+\gamma}{q}\Big)-\frac{r+\delta}{q}\Big|<\frac{\kappa c'}{q^{1+j}}\Big\}.$$
Furthermore,  for each $v=(p,r,q)\in\ZZ^2\times\NN$, we associate the interval
$$\Delta_{\bm\theta}(v):=\Big\{x\in I:\Big|x-\frac{p+\gamma}{q}\Big|<\frac{c'}{q^{1+i}}\Big\}.$$
Then,  with $\cA_0 \subset \cB_0 $ as in  \S\ref{winstrat}, the following is the inhomogeneous analogue of Lemma~\ref{L:Bad_c}.

\begin{lemma}\label{L:Bad_c'}
Let $\cA_0$, $  \sV $ and $ \Delta_{\bm\theta}(v) $ be as above.  Then
$$\cA_0\setminus\bigcup_{v\in\sV}\Delta_{\bm\theta}(v) \ \subset \ \Bad_{\bm\theta}^f(i,j)  \, . $$
\end{lemma}

\noindent{\em Proof. \ } The proof is similar to the homogeneous proof but is included for the sake of completeness.
Let $x\in\cA_0$. Suppose $x\notin\Bad^f_{\bm\theta}(i,j)$. Then there exists $v=(p,r,q)\in\ZZ^2\times\NN$ such that $$\Big|x-\frac{p+\gamma}{q}\Big| \ < \ \frac{c'}{q^{1+i}}, \qquad \Big|f(x)-\frac{r+\delta}{q}\Big| \ < \ \frac{c'}{q^{1+j}}.$$ In view of the fact that
$$\Big|x-\frac{p+\gamma}{q}\Big| \ < \ \frac{c'}{q^{1+i}} \ \le \  c \ \le \ \frac{l}{2},$$
it follows that $\frac{p+\gamma}{q}\in\cB_0\subset I$. Hence
\begin{eqnarray*}
\Big|f\Big(\frac{p+\gamma}{q}\Big)-\frac{r+\delta}{q}\Big| & \le & \Big|f\Big(\frac{p+\gamma}{q}\Big)-f(x)\Big|
+\Big|f(x)-\frac{r+\delta}{q}\Big|\\[1ex]
& \le & (\kappa-1)\Big|x-\frac{p+\gamma}{q}\Big|+\frac{c'}{q^{1+j}} \\[1ex]
& < & \frac{(\kappa-1)c'}{q^{1+i}}+\frac{c'}{q^{1+j}}  \ \le  \ \frac{\kappa c'}{q^{1+j}}.
\end{eqnarray*}
Thus $x\in\Delta_{\bm\theta}(v)$ and $v\in\sV$.  This completes the proof of the lemma.
\vspace*{-2ex}

\hfill $ \boxtimes $

\vspace*{2ex}

For $n\ge1$, let $$H'_n:=2c'l^{-1}R^n$$ and
$$\sV_n:=\{(p,r,q)\in\sV: H'_n\le q^{1+i}< H'_{n+1}\}.$$
Observe that  $H'_1=2c'l^{-1}R\le1$ and so it follows that $\sV=\bigcup_{n=1}^\infty\sV_n$. We inductively define a subtree $\cQ$ of $\cS'$ as follows.
Let $\cQ_0=\{\tau_0\}$. If $\cQ_{n-1}$ $(n\ge1)$ is defined, we let
$$\cQ_n:=\Big\{\tau\in\cS'_{\suc}(\cQ_{n-1}):\cI(\tau)\cap\bigcup_{v\in\sV_n}\Delta_{\bm\theta}(v)=\emptyset\Big\}.$$
Then $$\cQ:=\bigcup_{n=0}^\infty\cQ_n$$ is a subtree of $\cS'$ and, by
construction, we have that
\begin{equation*}\label{IE:tau}
\cI(\tau) \ \subset \ \cA_0\setminus\bigcup_{v\in\sV_n}\Delta_{\bm\theta}(v)   \qquad \forall  \ n\ge1 \quad {\rm and } \quad  \tau\in\cQ_n.
\end{equation*}

Armed with the following result, the same arguments as in \S\ref{howzthat} with the most obvious modifications enables us to prove Theorem \ref{T:inhomogeneous}.  In view of this the details of the proof of Theorem \ref{T:inhomogeneous} modulo Proposition \ref{keyinhom} are omitted.

\begin{proposition} \label{keyinhom}
The tree $\cQ$ has an $([R]-12)$-regular subtree.
\end{proposition}

In order to establish  the proposition, it suffices to prove the following statement.

\begin{lemma}\label{L:inhomo-curve}
For any $n\ge1$ and $\tau\in\cQ_{n-1}$, there is at most one $v\in\sV_n$ such that $\cI(\tau)\cap\Delta_{\bm\theta}(v)\ne\emptyset$. Moreover, $\rho(\Delta_{\bm\theta}(v))\le lR^{-n}$. Therefore,
$$\#\Big\{\tau'\in\cS'_{\suc}(\tau):\cI(\tau')\cap\bigcup_{v\in\sV_n}\Delta_{\bm\theta}(v)\ne\emptyset\Big\} \ \le \ 2.$$
\end{lemma}

\noindent{\em Proof. \ }
Suppose $v_s=(p_s,r_s,q_s)\in\sV_n$ and $\cI(\tau)\cap\Delta_{\bm\theta}(v_s)\ne\emptyset$, $s=1,2$. We need to prove that $v_1=v_2$. Without loss of generality, assume that $q_1\ge q_2$. Let $x_s\in\cI(\tau)\cap\Delta_{\bm\theta}(v_s)$. Then
$$\Big|x_s-\frac{p_s+\gamma}{q_s}\Big|<\frac{c'}{q_s^{1+i}} \ ,  \qquad s=1,2$$
and $$|x_1-x_2| \ \le \ \rho(\cI(\tau)) \ = \ lR^{-n+1}.$$
It follows that
\begin{eqnarray}
|(q_1-q_2)x_1-(p_1-p_2)| & = & \Big|q_1\Big(x_1-\frac{p_1+\gamma}{q_1}\Big)
-q_2\Big(x_2-\frac{p_2+\gamma}{q_2}\Big)-q_2(x_1-x_2)\Big|  \nonumber \\[1ex]
& \le & q_1\Big|x_1-\frac{p_1+\gamma}{q_1}\Big|+q_2\Big|x_2-\frac{p_2+\gamma}{q_2}\Big|+q_2|x_1-x_2| \nonumber  \\[1ex]
& \le & \frac{c'}{q_1^i}+\frac{c'}{q_2^i}+q_2lR^{-n+1} \nonumber  \\[1ex]
& \le & \frac{2c'}{q_2^i}+q_2lR^{-n+1}.   \label{w1}
\end{eqnarray}
Moreover,

\begin{eqnarray}
\Big| (q_1-q_2)f\Big(\frac{p_1+\gamma}{q_1}\Big)- (r_1-r_2) \Big| & = &
\Big| q_1\Big(f\Big(\frac{p_1+\gamma}{q_1}\Big)-\frac{r_1+\delta}{q_1}\Big)
 - q_2\Big(f\Big(\frac{p_2+\gamma}{q_2}\Big)-\frac{r_2+\delta}{q_2}\Big) \nonumber \\[1ex]
& &  \hspace*{6ex} -q_2\Big(f\Big(\frac{p_1+\gamma}{q_1}\Big)-f\Big(\frac{p_2+\gamma}{q_2}\Big)\Big) \Big| \nonumber \\[1ex]
& \le & q_1\Big|f\Big(\frac{p_1+\gamma}{q_1}\Big)-\frac{r_1+\delta}{q_1}\Big|+
q_2\Big|f\Big(\frac{p_2+\gamma}{q_2}\Big)-\frac{r_2+\delta}{q_2}\Big| \nonumber \\[1ex]
& &  \hspace*{6ex} + q_2\Big|f\Big(\frac{p_1+\gamma}{q_1}\Big)-f\Big(\frac{p_2+\gamma}{q_2}\Big)\Big| \nonumber  \\[1ex]
& \le & \frac{\kappa c'}{q_1^j}+\frac{\kappa c'}{q_2^j}+q_2\kappa\Big|\frac{p_1+\gamma}{q_1}-\frac{p_2+\gamma}{q_2}\Big| \nonumber  \\[1ex]
& \le & \frac{\kappa c'}{q_1^j}+\frac{\kappa c'}{q_2^j}+q_2\kappa\Big(\frac{c'}{q_1^{1+i}}+\frac{c'}{q_2^{1+i}}+lR^{-n+1}\Big) \nonumber \\[1ex]
& \le & \frac{4\kappa c'}{q_2^j}+q_2\kappa lR^{-n+1}.  \label{w2}
\end{eqnarray}

Suppose for the moment  that $q_1>q_2$  and let
$$P_0:=\Big(\frac{p_1-p_2}{q_1-q_2},\frac{r_1-r_2}{q_1-q_2}\Big)  \, . $$  We show that
\begin{equation}\label{vb206}
\cI(\tau)\cap\Delta(P_0)\ne\emptyset
\end{equation}
and that $P_0\in\cP$, where $\cP$ and $\Delta(P_0)$ are defined by \eqref{E:cP} and \eqref{E:Delta}, respectively.     In view of \eqref{w1},   it follows that
\begin{eqnarray*}
(q_1-q_2)^{1+i}\Big|x_1-\frac{p_1-p_2}{q_1-q_2}\Big| & \le   &  (q_1-q_2)^i \ \Big(\frac{2c'}{q_2^i}+q_2lR^{-n+1}\Big)  \\[1ex] &\le&  2c'\frac{q_1^i}{q_2^i}+q_1^{1+i}lR^{-n+1}   \ \le \   2c'R+2c'R^2  \ <  \ c.
\end{eqnarray*}
So $x_1\in\Delta(P_0)$ and \eqref{vb206} follows.
 Also, the above inequality implies that
$$\Big|x_1-\frac{p_1-p_2}{q_1-q_2}\Big| \ < \ c \ \le \ \frac{l}{2}  \, $$
and since $x_1\in\cA_0$, it follows that $\frac{p_1-p_2}{q_1-q_2}\in\cB_0\subset I$. Moreover, by making use
of \eqref{w2} we have that
\begin{align*}
(q_1-q_2)^{1+j} \Big|f\Big(\frac{p_1-p_2}{q_1-q_2}\Big) & -\frac{r_1-r_2}{q_1-q_2}\Big|
\ \le  \ (q_1-q_2)^{1+j}\Big|f\Big(\frac{p_1-p_2}{q_1-q_2}\Big)-f\Big(\frac{p_1+\gamma}{q_1}\Big)\Big| \\[1ex]
&~ \hspace{19ex}  + (q_1-q_2)^j\Big|(q_1-q_2)f\Big(\frac{p_1+\gamma}{q_1}\Big)-(r_1-r_2)\Big|\\[2ex]
\le \  &   (q_1-q_2)^{1+j}\kappa\Big|\frac{p_1-p_2}{q_1-q_2}-\frac{p_1+\gamma}{q_1}\Big|\\[1ex]
&  \hspace{9ex} + (q_1-q_2)^j\Big(\frac{4\kappa c'}{q_2^j}+q_2\kappa lR^{-n+1}\Big)\\[1ex]
=   \ & (q_1-q_2)^j\Big(q_2\kappa\Big|\frac{p_1+\gamma}{q_1}-\frac{p_2+\gamma}{q_2}\Big|+\frac{4\kappa c'}{q_2^j}+q_2\kappa lR^{-n+1}\Big)\\[1ex]
\le \ & q_1^j\Big(q_2\kappa\Big(\frac{c'}{q_1^{1+i}}+lR^{-n+1}+\frac{c'}{q_2^{1+i}}\Big)+\frac{4\kappa c'}{q_2^j}+q_2\kappa lR^{-n+1}\Big)\\[1ex]
\le \ & 6\kappa c'\frac{q_1^j}{q_2^j}+2q_1^{1+j}\kappa lR^{-n+1}\\[1ex]
\le \ & 6\kappa c'R+4\kappa c'R^2 \ < \ \kappa c.
\end{align*}
Thus $P_0\in\cP$ and so there exists a  unique integer $n_0\ge1$  such that $P_0\in\cP_{n_0}$.   Suppose for the moment that  $n_0\le n-1$.  Then there exists $\tau'\in\cS_{n_0}$ such that $\tau\prec\tau'$  and it follows that $\cI(\tau)\cap\Delta(P_0) \ \subset  \  \cI(\tau')\cap\Delta(P_0)=\emptyset$ contrary to \eqref{vb206}.
Thus, $n_0\ge n $ and so
$$H(P_0) \ \ge \  H_{n_0} \ \ge \  H_{n} \ = \  42\kappa^3cl^{-1}R^n. $$
On the other hand, we have that
$$H(P_0) \ \le \  \kappa(q_1-q_2)^{1+i} \ \le \  \kappa q_1^{1+i} \ \le \  \kappa H'_{n+1} \ = \ \frac{1}{5}\kappa cl^{-1}R^{n-1}  \, . $$
This contradicts the above lower bound for $H(P_0)$ and we conclude that $q_1=q_2$.  Since $q_2\le H'_{n+1}\le cl^{-1}R^{n-1}$, it now  follows from \eqref{w1} and \eqref{w2} that
$$|p_1-p_2| \ \le \  \frac{2c'}{q_2^i}+q_2lR^{-n+1} \ \le \  2c'+c \ < \ 1$$ $$$$
and
$$|r_1-r_2| \ \le \  \frac{4\kappa c'}{q_2^j}+q_2\kappa lR^{-n+1} \ \le \  4\kappa c'+c\kappa \ < \ 1.$$
Thus,  $p_1=p_2$ and $r_1=r_2$. In other words,  $v_1=v_2$ and this proves the  main  substance  of the lemma.  To prove the `moreover' part, it is easily verified that for any $v\in\sV_n$  we have  that
$$\rho(\Delta_{\bm\theta}(v)) \ = \ \frac{2c'}{q^{1+i}} \ \le  \ \frac{2c'}{H'_n} \ = \ lR^{-n}.$$
The `therefore' part of the lemma is a direct consequence of  this and the  fact that there is at most one $v\in\sV_n$ such that $\cI(\tau)\cap\Delta_{\bm\theta}(v)$ is non-empty.
\vspace*{-2ex}

\hfill $ \boxtimes $

\vspace*{2ex}

As already mentioned, given Proposition \ref{keyinhom},   the proof of Theorem \ref{T:inhomogeneous} follows on adapting the arguments of \S\ref{howzthat}.

\section{The proof of  Theorem \ref{thm2} \label{thm2strat}}

The basic strategy towards establishing the winning result for lines is the same as when considering curves. To begin with observe that for any line  $\L_{a,b}$ given by
$$y=f(x) := ax +b $$ and
 $\bm\theta \in \RR^2$,
$$
\Bad_{\bm\theta}^f(i,j)  \, :=  \, \{x\in \RR:(x,f(x))\in\Bad_{\bm\theta}(i,j)\}  \; =  \; \pi( \Bad_{\bm\theta}(i,j) \cap \L_{a,b} ) \, .
$$
As in the case of curves, without loss of generality we will assume that $j \leq i$.   Thus,  the
homogeneous  case of Theorem \ref{thm2} is easily deduced from the following  statement.

\begin{theorem}\label{T2homo:main}
Let  $(i,j)$ be a pair of real numbers
satisfying   $0<j\le i<1$ and $i+j=1$.   Given $a,b\in \RR$, suppose there exists
$\epsilon>0$ such that
\begin{equation} \label{diocondbetterT2}
\liminf_{q \to \infty}   q^{\frac{1}{j}-\epsilon}  \max \{\|q a
\|, \|q b
\|  \}  > 0 \, .
\end{equation}
Then $\Bad^f(i,j)$
is a $1/2$-winning subset of $\RR$. Moreover, if $ a \in \QQ$ the statement is true with $ \epsilon=0$ in  (\ref{diocondbetterT2}).
\end{theorem}

\noindent Note that in view of Remark 5  after the statement of Theorem \ref{thm2},    we do not require that $ a \neq 0$ in Theorem \ref{T2homo:main}   since   $ j \le i $.

\subsection{The winning strategy for Theorem \ref{T2homo:main} \label{thm2strat2}}   Let $\beta\in(0,1)$. We want to prove that $\Bad^f(i,j)$ is $(\frac{1}{2},\beta)$-winning. In the first round of the game, \textbf{B}hupen chooses a closed interval $\cB_0\subset \RR$.  Now  \textbf{A}yesha chooses any  closed interval $\cA_0\subset\cB_0$ with diameter $\rho(\cA_0)=\frac{1}{2} \rho(\cB_0)$.
Let
$$R:=( 2 \beta^{-1})^4 \, ,   \qquad l:=  \rho(\cA_0)  \qquad {\rm and }  \qquad c_1:=c_0\max\{|x|_{\max}+l,1\}^{-1}  \,  $$
where
$$
|x|_{\max} := \max\{ |x| : x \in \cA_0 \} \qquad {\rm and }  \qquad  c_0:=\inf_{q \in \NN}   q^{\frac{1}{j}-\epsilon}  \max \{\|q a
\|, \|q b
\|  \}  > 0  \, .
 $$
The fact that $c_0 > 0$ follows from the Diophantine condition \eqref{diocondbetterT2}.  Recall, that by hypothesis  $\epsilon=0$ if $a$ is rational and $\epsilon>0$ otherwise. We denote
$$\kappa:=|a|+1,$$
and choose   $\mu\ge1$  such that
\begin{equation}\label{E:mulines}
20\kappa^2R^{\frac{3}{i}-i\mu}<1.
\end{equation}
If $a\in\QQ$ , we also require that
\begin{equation}\label{E:mulines2}
R^{\mu-1}\ge\kappa d^2
\end{equation}
where  $d\in\NN$ is the smallest positive integer such that $da\in\ZZ$. Next, if $a\notin\QQ$, so that $\epsilon>0$, we let $\lambda>0$ be such that
\begin{equation}\label{E:lambdalines}
R^{\lambda j^{-n}} \, \ge  \, \kappa c_1^{-1}R^{\frac{1+i}{j\epsilon}n} \quad {\rm for }  \quad n \ge1.
\end{equation}
If $a\in\QQ$, we simply let $\lambda=0$.
In turn, let
$$\lambda_0:=0   \qquad  {\rm and }  \qquad \lambda_k  \, := \, \lambda \, j^{-k}+\frac{k}{i}+\mu  \quad {\rm for }  \quad  k\ge1,$$
and let
\begin{equation} \label{E:clines}
c:=\min\Big\{\frac{c_1}{4\kappa}lR^{-1}, l^{-i}, \frac{1}{8\kappa}R^{-2-\frac{\lambda_1}{1+i}}\Big\} \, .
\end{equation}

\vspace*{3ex}

For each rational point $ P=(\frac{p}{q},\frac{r}{q}) \in \RR^2$ we associate the interval $$\Delta(P):=\Big\{x\in\RR:\Big|x-\frac{p}{q}\Big|<\frac{c}{q^{1+i}}\Big\} $$
and  we let
$$\cP:=\Big\{P=\Big(\frac{p}{q},\frac{r}{q}\Big)  \,  :  \, \cA_0\cap\Delta(P) \neq \emptyset, \,   \Big|b+\frac{ap-r}{q}\Big|<\frac{\kappa c}{q^{1+j}}\Big\}.$$

The following inclusion is a simple consequence of the manner in which the above quantities  and   objects have been defined.

\begin{lemma}\label{L:Bad_clines}  Let $\cA_0$, $  \cP  $ and $ \Delta(P) $ be as above.  Then
$$\cA_0\setminus\bigcup_{P\in\cP}\Delta(P) \ \subset \ \Bad^f(i,j)  \, . $$
\end{lemma}

\noindent{\em Proof. \ }
Let $x\in\cA_0$. Suppose $x\notin\Bad^f(i,j)$. Then there exists $P=(\frac{p}{q},\frac{r}{q})\in\QQ^2$ such that $$\Big|x-\frac{p}{q}\Big| \ < \ \frac{c}{q^{1+i}}, \qquad \Big|f(x)-\frac{r}{q}\Big| \  =  \ \Big|ax+b-\frac{r}{q}\Big| \  < \ \frac{c}{q^{1+j}}.$$ It follows that

\begin{eqnarray*}
\Big|b+\frac{ap-r}{q}\Big| & =  & \Big|\Big(ax+b-\frac{r}{q}\Big)-a\Big(x-\frac{p}{q}\Big)  \Big|  \\[1ex]
& \le  &   \frac{c}{q^{1+j}}+|a|\frac{c}{q^{1+i}}  \\[1ex]
& \le   &  \frac{\kappa c}{q^{1+j}}     \, .
\end{eqnarray*}
The upshot is that $x\in\Delta(P)$ with $P\in\cP$. This completes the proof of the lemma.

\vspace*{-2ex}

\hfill $ \boxtimes $

\vspace*{2ex}

Next, just as in \S\ref{winstrat}, but with $\cA_0$ and $\cP$ as above, let
\begin{itemize}
  \item  $\cT$ be an $[R]$-regular rooted tree with root $\tau_0$ ,
  \item $\cI$ be an injective map from $\cT$ to the set of closed subintervals of $\cA_0$ ,
  \item $\cS=\bigcup_{n=0}^\infty\cS_n$ be  a subtree of $\cT$ associated with a partition $\cP_n$ of $\cP$.
\end{itemize}
The following proposition is the lines analogue of Proposition \ref{P:main}.  It enables us to deduce Theorem \ref{T2homo:main} (and thus the homogeneous case of Theorem \ref{thm2})   by  adapting the arguments  of \S\ref{howzthat} in the most obvious manner. In view of this the details of the proof of Theorem~\ref{T2homo:main} modulo Proposition \ref{P:mainlines} are omitted.

\begin{proposition}\label{P:mainlines}
There exists a partition $\cP=\bigcup_{n=1}^\infty\cP_n$ such that the  tree $\cS$ has an $([R]-5)$-regular subtree.
\end{proposition}

\subsubsection{Proof of Proposition \ref{P:mainlines}}
 As in the `curves' proof,  to each point $P=(\frac{p}{q},\frac{r}{q})\in\cP$, we attach a rational line
$$L_P:=\{(x,y)\in\RR^2:Ax+By+C=0\}$$
passing through $P$ where  $A,B,C\in\ZZ$  are co-prime  with  $(A,B)\ne(0,0)$ and such that $$|A|\le q^i  \,   \quad {\rm and }  \quad  |B|\le q^j  \, . $$
Associated with each point $P\in\cP$, we also consider the quantity
$$E:=A+Ba  \, . $$
Then
\begin{equation} \label{iona1}
|E| \ \le \  q^i+|a|q^j\le\kappa q^i   \, .
\end{equation}
Note  that if $x\in\Delta(P)$, then
\begin{equation} \label{iona2}
|Ex+Bb+C| \ = \   \Big|E\Big(x-\frac{p}{q}\Big)+B\Big(b+\frac{ap-r}{q}\Big)\Big| \ < \ \frac{2\kappa c}{q}   \, .
\end{equation}

The following statement enables us to construct the desired partition in Proposition~\ref{P:mainlines}.

\begin{lemma}\label{L:qEbound}
For any $P=(\frac{p}{q},\frac{r}{q})  \in \cP$, we have
\begin{equation}\label{iona3}
q^{1-j\epsilon}|E|\ge c_1.
\end{equation}
\end{lemma}

\noindent{\em Proof. \ }  If $B=0$, then $q^{1-j\epsilon}|E|=q^{1-j\epsilon}|A|\ge1\ge c_1$ and we are done. If $B\ne0$, then it is easily verified that $$q^{1-j\epsilon}\max\{|E|,|Bb+C|\} \ \ge \ |B|^{\frac{1}{j}-\epsilon}\max\{|Ba+A|,|Bb+C|\} \ \ge \  c_0.$$
If $|E|$ is the maximum in the above then again we are done.  So suppose $q^{1-j\epsilon}|Bb+C|\ge c_0$. Since $P\in\cP$, there exists $x\in\cA_0\cap\Delta(P)$ and it follows that
 \begin{eqnarray*}
 q^{1-j\epsilon}|Ex|  & \ge   &   q^{1-j\epsilon}|Bb+C|-q^{1-j\epsilon}|Ex+Bb+C| \\[1ex]
 & \stackrel{\eqref{iona2}}{\ge}  &    c_0-2\kappa c \ \ge  \  c_1(|x|_{\max}+l)-c_1l  \ =  \ c_1|x|_{\max}.
\end{eqnarray*}   This proves the lemma.

\vspace*{-2ex}

\hfill $ \boxtimes $

\vspace*{2ex}

A particular consequence of \eqref{iona3} is that $E\ne0$. Thus every line $L_P$ intersects at the line $\L_{a,b}$ given by $y=f(x)= ax +b $  at a single point.

For $n\ge1$, let $$H_n:=4\kappa cl^{-1}R^n$$
and
$$\cP_n:=\Big\{P=\Big(\frac{p}{q},\frac{r}{q}\Big)\in\cP: H_n\le q|E|< H_{n+1}\Big\}.$$
Note that if $P\in\cP_n$, then
\begin{equation}\label{iona4}\kappa q^{1+i}\stackrel{\eqref{iona1}}{\ge} q|E|\ge H_n \, .\end{equation}
Next let
$$\cP_{n,k}:=\{P\in\cP_n: H_nR^{\lambda_{k-1}}\le \kappa q^{1+i}< H_nR^{\lambda_k}\}   \quad {\rm for }  \quad 1\le k\le n.$$

\begin{lemma}\label{L:partitionlines}  With $\cP_n$ and $\cP_{n,k}$ as above, we have that
$$\cP=\bigcup_{n=1}^\infty\cP_n  \quad  { and }  \quad  \cP_n=\bigcup_{k=1}^{n}\cP_{n,k} \, .  $$
\end{lemma}

\noindent{\em Proof. \ }  Note that by \eqref{iona3}, for any $P \in \cP$  we have that   $q|E|\ge c_1$ and by definition
\begin{equation}\label{E:H_1}
H_1=4\kappa cl^{-1}R\le c_1.
\end{equation}
Thus, $\cP=\bigcup_{n=1}^\infty\cP_n$.

To prove the second conclusion, we first show that
\begin{equation}\label{E:second-partition}
\kappa q^{1+i} < H_nR^{\lambda_n} \quad {\rm for }  \quad P \in \cP_n.
\end{equation}
We consider two separate cases. If $a\notin\QQ$, on combining the fact that $ q|E|<H_{n+1}$ with \eqref{iona3} implies that
$$
q^{j\epsilon}<c_1^{-1}H_{n+1}.
$$
It then follows that
$$\frac{\kappa q^{1+i}}{H_nR^{\lambda_n}} \ < \ \frac{\kappa(c_1^{-1}H_{n+1})^{\frac{1+i}{j\epsilon}}}{H_nR^{\lambda_n}}
 \ = \ \frac{\kappa c_1^{-1}(c_1^{-1}H_1R^n)^{\frac{1+i}{j\epsilon}-1}}{R^{\lambda_n-1}}
 \ \stackrel{\eqref{E:H_1}}{\le} \ \kappa c_1^{-1}R^{(\frac{1+i}{j\epsilon}-1)n-\lambda j^{-n}} \ \stackrel{\eqref{E:lambdalines}}{\le} \ 1.$$
If $a\in\QQ$, then $|E|\ge 1/d$, and hence the fact $ q|E|<H_{n+1}$ implies that $q<dH_{n+1}$. It follows that
$$\frac{\kappa q^{1+i}}{H_nR^{\lambda_n}} \ < \ \frac{\kappa d^2H_{n+1}^{1+i}}{H_nR^{\frac{n}{i}+\mu}} \ \stackrel{\eqref{E:mulines2}}{\le} \
H_1^{i}R^{n(i-\frac{1}{i})} \ < \ 1.$$
This proves \eqref{E:second-partition}. Now \eqref{E:second-partition} together with \eqref{iona4} implies that  $\cP_n=\bigcup_{k=1}^{n}\cP_{n,k}$.
~\vspace*{-1ex}
 \hfill $ \boxtimes $

\vspace*{3ex}

We claim that the partition of $\cP$  given by Lemma \ref{L:partitionlines} satisfies the requirement of Proposition \ref{P:mainlines}.   The key towards establishing the claim is the following lemma.  It is the lines analogue of Lemma \ref{L:const}.

\begin{lemma}\label{L:constlines}
For any $n\ge1$, $1\le k\le n$ and $\tau\in\cS_{n-k}$, the map $P\mt L_P$ is constant on $$\cP_{n,k}(\tau):=\{P\in\cP_{n,k}:\cI(\tau)\cap\Delta(P)\ne\emptyset\}.$$
\end{lemma}

\noindent{\em Proof. \ }
Let $P_1=(\frac{p_1}{q_1},\frac{r_1}{q_1})$ and $P_2=(\frac{p_2}{q_2},\frac{r_2}{q_2})$ be distinct points in $\cP_{n,k}(\tau)$. We need to prove that $L_{P_1}=L_{P_2}$.  We let $A_s$, $B_s$, $C_s$ and $E_s$ be the respective quantities associated with $P_s$, $s=1,2$, and consider two separate cases.
\vspace*{2ex}

{\em Case (1). \ } Suppose $k=1$. Then

\begin{eqnarray*}
|A_1p_2+B_1r_2+C_1q_2| & =&q_2\Big|A_1\Big(\frac{p_2}{q_2}-\frac{p_1}{q_1}\Big)+B_1\Big(\frac{r_2}{q_2}-\frac{r_1}{q_1}\Big)\Big| \\[1ex]
&=&q_2\Big|E_1\Big(\frac{p_2}{q_2}-\frac{p_1}{q_1}\Big)+B_1\Big(\frac{ap_1-r_1}{q_1}-\frac{ap_2-r_2}{q_2}\Big)\Big|\\[1ex]
&\le& q_2|E_1|\Big(\frac{c}{q_1^{1+i}}+\frac{c}{q_2^{1+i}} +lR^{-n+1}\Big) \\[1ex]
 & ~ & ~ \hspace{17ex} + \ q_2|B_1|\Big(\frac{\kappa c}{q_1^{1+j}}+
\frac{\kappa c}{q_2^{1+j}}\Big)\\[1ex]
&\le& \frac{q_2}{q_1}q_1|E_1|\Big(\frac{2\kappa c}{H_n}+lR^{-n+1}\Big)+\kappa c \, \Big(\frac{q_2}{q_1}+\frac{q_1^j}{q_2^j}\Big)\\[1ex]
&\le & R^{\frac{\lambda_1}{1+i}}H_{n+1}\cdot\frac{2\kappa c(1+2R)}{H_n}+2\kappa cR^{\frac{\lambda_1}{1+i}}\\[1ex]
&< & 8\kappa cR^{2+\frac{\lambda_1}{1+i}} \ \le \ 1.
\end{eqnarray*}
Since the left hand side of the above inequality is an integer, it follows that
$A_1p_2+B_1r_2+C_1q_2=0$. Similarly, we obtain that $A_2p_1+B_2r_1+C_2q_1=0$.   The upshot is that both the lines $ L_{P_1} $and $L_{P_2} $  pass through both the points $P_1$ and $P_2$,  and so we must have that  $L_{P_1}=L_{P_2}$.

\vspace*{2ex}

{\em Case (2). \ } Suppose $k\ge2$.    We prove that  $L_{P_1} = L_{P_2}$ by contradiction.   Thus, assume that $L_{P_1}\ne L_{P_2}$.   We first establish various preliminary estimates.  Let
$$m_q:=(\kappa^{-1}H_nR^{\lambda_{k-1}})^{\frac{1}{1+i}} \qquad {\rm and}   \qquad  M_q:=(\kappa^{-1}H_nR^{\lambda_k})^{\frac{1}{1+i}}.$$
Then, by definition,  for $P=(\frac{p}{q},\frac{r}{q})\in\cP_{n,k}$ we have that  $m_q\le q<M_q$.  Also let
$$M_E:=m_q^{-1}H_{n+1}  \qquad {\rm and}   \qquad     M_B:=M_q^j.$$
Then $|E_s|\le M_E$, $|B_s|\le M_B$, $s=1,2$. We claim that
\begin{equation} \label{iona6}
M_E^{1+\frac{1}{i}} \ \le \  M_B^jM_E^{2+j} \ \le \  M_B^{1+i}M_E^{1+i} \ < \ \frac{1}{5}cl^{-1}R^{n-k} \, .
\end{equation}
First observe that
$$\frac{M_E}{M_B^{i/j}} \ = \ \frac{H_{n+1}}{m_qM_q^i} \ = \ \kappa R^{1-\frac{\lambda_{k-1}+i\lambda_k}{1+i}} \ \le \ \kappa R^{1-\mu}\le1.$$
It then follows that
$$\frac{M_E^{1+\frac{1}{i}}}{M_B^jM_E^{2+j}} \ = \ \Big(\frac{M_E}{M_B^{i/j}}\Big)^{\frac{j^2}{i}} \ \le \ 1   \qquad  {\rm and}  \qquad
\frac{M_B^jM_E^{2+j}}{M_B^{1+i}M_E^{1+i}} \ = \ \Big(\frac{M_E}{M_B^{i/j}}\Big)^{2j} \ \le \ 1.$$
This establishes the left  and middle  inequalities within  \eqref{iona6}.  Regarding the right inequality, we have that
\begin{eqnarray*}
M_B^{1+i}M_E^{1+i} & = & \Big(\frac{H_{n+1}}{M_q^i}\cdot\frac{M_q}{m_q}\Big)^{1+i}=
\frac{\kappa^iH_{n+1}^{1+i}}{(H_nR^{\lambda_k})^i}\cdot R^{\lambda_k-\lambda_{k-1}} \\[1ex]
& = & \kappa^iH_nR^{j\lambda_k-\lambda_{k-1}+1+i} \ = \ 4\kappa^{1+i}cl^{-1}R^{n-k+\frac{1}{i}-i\mu+1+i}\\[1ex]
 & <& \frac{1}{5}cl^{-1}R^{n-k}.
\end{eqnarray*}
With $s=1$ or $2$, let $(x_s,ax_s+b)$ denote the  intersection point  of $L_{P_s}$ with the line $\L_{a,b}$.
Then
$$E_s\Big(x_s-\frac{p_s}{q_s}\Big)+B_s\Big(b+\frac{ap_s-r_s}{q_s}\Big) \ = \ 0  \, , $$
and so
$$\Big|x_s-\frac{p_s}{q_s}\Big| \ = \ \frac{|B_s|}{|E_s|}\Big|b+\frac{ap_s-r_s}{q_s}\Big| \ \le \ \frac{|B_s|}{|E_s|}\frac{\kappa c}{q_s^{1+j}}, \qquad  \quad  s=1,2.$$
Hence
\begin{eqnarray} \label{iona7} |x_1-x_2| & \le & \Big|x_1-\frac{p_1}{q_1}\Big|+\Big|x_2-\frac{p_2}{q_2}\Big|+\Big|\frac{p_1}{q_1}-\frac{p_2}{q_2}\Big| \nonumber \\[1ex] \nonumber
& \le & \frac{|B_1|}{|E_1|}\frac{\kappa c}{q_1^{1+j}}+\frac{|B_2|}{|E_2|}\frac{\kappa c}{q_2^{1+j}}+
\Big(\frac{c}{q_1^{1+i}}+\frac{c}{q_2^{1+i}}+lR^{-n+k}\Big)  \\[1ex]  \nonumber
& \le  & \frac{\kappa c}{q_1|E_1|}+\frac{\kappa c}{q_2|E_2|}+\Big(\frac{\kappa c}{q_1|E_1|}+\frac{\kappa c}{q_2|E_2|}+lR^{-n+k}\Big)\\[1ex]  \nonumber
& \le   & \frac{4\kappa c}{H_n}+lR^{-n+k} \ = \ lR^{-n}+lR^{-n+k}  \\[1ex]
& \le   & 2lR^{-n+k}  \, .
\end{eqnarray}
This completes the preliminaries.  Recall that we are assuming that  $L_{P_1}\ne L_{P_2}$ and the name of the game is to obtain a contradiction.  We first consider the case that
$L_{P_1}$ is parallel to $L_{P_2}$. Then there exist $(A,B)\in\ZZ^2\setminus\{(0,0)\}$ and nonzero integers $t_1,t_2$ such that $$(A_1,B_1)=t_1(A,B) \qquad {\rm and }  \qquad (A_2,B_2)=t_2(A,B).$$
Thus $$x_s \ = \ -\frac{B_sb+C_s}{B_sa+A_s} \ = \ -\frac{1}{Ba+A}\Big(Bb+\frac{C_s}{t_s}\Big)  \qquad  \quad  s=1,2 $$
and it follows that
 $$|x_1-x_2| \ = \ \frac{|t_1C_2-t_2C_1|}{|t_1t_2||Ba+A|} \ \ge \ \frac{1}{|t_1t_2||Ba+A|} \ = \ \frac{1}{|t_1E_2|}.$$
This together with \eqref{iona7} implies that
\begin{equation}  \label{iona8} |t_1E_2|\ge\frac{1}{2}l^{-1}R^{n-k}.\end{equation}
If $B_1=0$, then $|t_1|\le|A_1|\le|A_1|^{1/i}=|E_1|^{1/i}$. So
$$|t_1E_2| \ \le  \  M_E^{1+\frac{1}{i}} \ \stackrel{\eqref{iona6}}{<}  \ \frac{1}{5}cl^{-1}R^{n-k}$$
and this contradicts \eqref{iona8}. If $B_1\ne0$, then $|t_1|\le|B_1|$ and
$$|t_1E_2|^{1+i} \ \le \ (M_BM_E)^{1+i} \ \stackrel{\eqref{iona6}}{<}  \ \frac{1}{5}cl^{-1}R^{n-k}.$$
This together with \eqref{iona8}  implies that
$$\Big(\frac{1}{2}l^{-1}R^{n-k}\Big)^{1+i} < \ \frac{1}{5}cl^{-1}R^{n-k}.$$
However, this contradicts the fact that $c\le l^{-i}$.  Hence, if $L_{P_1}$ is parallel to $L_{P_2}$ then we must have that $L_{P_1}= L_{P_2}$.

Now suppose $L_{P_1}$ is not parallel to $L_{P_2}$. Let $P_0=(\frac{p_0}{q_0},\frac{r_0}{q_0})\in\QQ^2$ be the intersection of $L_{P_1}$ and $L_{P_2}$. Then it follows that $A_1B_2-A_2B_1  $ is a nonzero integer and is divisible by $q_0$ and so
 \begin{equation} \label{1svlines}  q_0 \, \le \, |A_1B_2-A_2B_1|  \, = \, | E_1B_2 -  E_2B_1  |  \, .\end{equation}
 We first prove   that  $\Delta(P_1)\subset\Delta(P_0)$  and  that $P_0\in\cP$.  In view of the fact that
$$E_s\Big(\frac{p_s}{q_s}-\frac{p_0}{q_0}\Big)+B_s\Big(\frac{r_s-ap_s}{q_s}-\frac{r_0-ap_0}{q_0}\Big) \ = \ 0 \qquad s=1,2 $$
it is easily verified that
$$-(E_1B_2 -  E_2B_1)\Big(\frac{p_1}{q_1}-\frac{p_0}{q_0}\Big) \ = \ B_1B_2\Big(\frac{r_1-ap_1}{q_1}-\frac{r_2-ap_2}{q_2}\Big)
+B_1E_2\Big(\frac{p_1}{q_1}-\frac{p_2}{q_2}\Big).$$
It then follows that
\begin{eqnarray*}q_0\Big|\frac{p_1}{q_1}-\frac{p_0}{q_0}\Big| & \stackrel{\eqref{1svlines}}{\le} & |B_1B_2|\Big(\frac{\kappa c}{q_1^{1+j}}+\frac{\kappa c}{q_2^{1+j}}\Big)
+|B_1E_2|\Big(\frac{c}{q_1^{1+i}}+\frac{c}{q_2^{1+i}}+lR^{-n+k}\Big)  \\[1ex]
&\le  & |B_2E_1|\frac{\kappa c}{q_1|E_1|}+|B_1E_2|\frac{\kappa c}{q_2|E_2|} \\[1ex]
 & ~ & ~ \hspace{20ex} +|B_1E_2|\Big(\frac{\kappa c}{q_1|E_1|}+\frac{\kappa c}{q_2|E_2|}+lR^{-n+k}\Big)  \\[1ex]
& \le  &  M_BM_E\Big(\frac{4\kappa c}{H_n}+lR^{-n+k}\Big)  \\[1ex]
& \le  & 2M_BM_ElR^{-n+k}.
\end{eqnarray*}
So if $x\in\Delta(P_1)$, then
\begin{eqnarray*}  q_0^{1+i}\Big|x-\frac{p_0}{q_0}\Big|  & \le  &  q_0^{1+i}\Big|x-\frac{p_1}{q_1}\Big|+q_0^{1+i}\Big|\frac{p_1}{q_1}-\frac{p_0}{q_0}\Big|\\[1ex]
& \le  &  q_0^{1+i}\frac{c}{q_1^{1+i}}+2q_0^iM_BM_ElR^{-n+k} \\[1ex]
 & \stackrel{\eqref{1svlines}}{\le}  &  4M_B^{1+i}M_E^{1+i}\frac{\kappa c}{H_n}+4M_B^{1+i}M_E^{1+i}lR^{-n+k} \\[1ex]
& \le  & 5M_B^{1+i}M_E^{1+i}lR^{-n+k} \, \stackrel{\eqref{iona6}}{<}  \, c.
\end{eqnarray*}
Thus $x\in\Delta(P_0)$  and the upshot is that $\Delta(P_1)\subset\Delta(P_0)$. In turn,
since $\cA_0\cap\Delta(P_1)\ne~\emptyset$,  it follows that  $\cA_0\cap\Delta(P_0)\ne\emptyset$.  In view of this, in order to prove that  $P_0\in\cP$ we need to show that
\begin{equation}  \label{iona9} \Big|b+\frac{ap_0-r_0}{q_0}\Big| \ < \ \frac{\kappa c}{q_0^{1+j}}.\end{equation}
Since
$$E_s\Big(x_s-\frac{p_0}{q_0}\Big)+B_s\Big(b+\frac{ap_0-r_0}{q_0}\Big) \ = \ 0   \qquad  \quad  s=1,2 $$
we have that  $$|E_1B_2 -  E_2B_1|\Big|b+\frac{ap_0-r_0}{q_0}\Big| \ = \ |E_1E_2||x_1-x_2| \ \stackrel{\eqref{iona7}}{\le} \ 2M_E^2lR^{-n+k}.$$
Hence
\begin{eqnarray*}q_0^{1+j}\Big|b+\frac{ap_0-r_0}{q_0}\Big| & \stackrel{\eqref{1svlines}}{\le}  & |E_1B_2 -  E_2B_1|^{1+j}\Big|b+\frac{ap_0-r_0}{q_0}\Big|  \\[1ex]
& \le &  4M_B^jM_E^{2+j}lR^{-n+k}   \\[1ex]
& \stackrel{\eqref{iona6}}{<}   & \kappa c  \,
\end{eqnarray*}
and this established \eqref{iona9}.    Thus $P_0\in\cP$ and so there exists a  unique integer $n_0\ge1$  such that $P_0\in\cP_{n_0}$.   Suppose for the moment that  $n_0\le n-k$.  Then there exists $\tau'\in\cS_{n_0}$ such that $\tau\prec\tau'$, and hence $$\cI(\tau)\cap\Delta(P_1) \ \subset  \  \cI(\tau')\cap\Delta(P_0)=\emptyset.$$ This contradicts the fact that $P_1\in\cP_{n,k}(\tau)$.  Thus
$$n_0\ge n-k+1 \, , $$
and so
 $$ q_0^{1+i} \ \stackrel{\eqref{iona4}}{\ge} \  \kappa^{-1}H_{n_0} \ \ge \  \kappa^{-1}H_{n-k+1} \ = \ 4cl^{-1}R^{n-k+1}   \, . $$
On the other hand, we have that
$$q_0^{1+i} \ \le \ 4(M_BM_E)^{1+i} \ \stackrel{\eqref{iona6}}{<} \ cl^{-1}R^{n-k}  \, . $$
This contradicts the above lower bound for $q_0^{1+i}$ and so completes the proof of Case (2) and indeed the lemma.

\vspace*{-2ex}

\hfill $ \boxtimes $

\vspace*{2ex}

An important consequence of Lemma \ref{L:constlines} is the following lines analogue of Corollary~\ref{C:est}.

\begin{corollary}\label{line-C:est}
For any $n\ge1$, $1\le k\le n$ and $\tau\in\cS_{n-k}$, we have
$$\#\Big\{\tau'\in\cT_n:\cI(\tau')  \ \cap\bigcup_{P\in\cP_{n,k}(\tau)}\Delta(P)\ne\emptyset\Big\} \ \le \  2.$$
\end{corollary}

\noindent{\em Proof. \ }
By Lemma \ref{L:constlines} and \eqref{iona2}, there exists $(A,B,C)\in\ZZ^3$ with $E:=A+Ba\ne0$ such that for any $P=(\frac{p}{q},\frac{r}{q})\in\cP_{n,k}(\tau)$ and $x\in\Delta(P)$,
$$|Ex+Bb+C|<\frac{2\kappa c}{q}   \qquad {\rm   and }  \qquad q|E|\ge H_n \, . $$
Thus
$$\Big|x+\frac{Bb+C}{E}\Big| \ < \ \frac{2\kappa c}{q|E|} \ \le \ \frac{2\kappa c}{H_n} \ = \ \frac{1}{2}lR^{-n}.$$
This implies that $\bigcup_{P\in\cP_{n,k}(\tau)}\Delta(P)$ is contained in the open interval
\begin{equation}\label{E:interval}
\Big(-\frac{Bb+C}{E}-\frac{1}{2}lR^{-n},-\frac{Bb+C}{E}+\frac{1}{2}lR^{-n}\Big),
\end{equation}
which has length $lR^{-n}$. Since the intervals $\{\cI(\tau'):\tau'\in\cT_n\}$ are of length $lR^{-n}$ and have mutually disjoint interiors, there can be at most $2$ of them that intersect the interval \eqref{E:interval}. This proves the corollary.
\vspace*{-2ex}

\hfill $ \boxtimes $

\vspace*{2ex}

We are now in the position to prove Proposition \ref{P:mainlines}.   In view of Proposition \ref{propJA},  it suffices to prove that the intersection of $\cS$ with every $6$-regular subtree of $\cT$ is infinite.     Let $\cR\subset\cT$ be a $6$-regular subtree and let
$$\cR':=\cR\cap\cS   \qquad {\rm   and }  \qquad a_n:=\#\cR'_n  \quad (n \ge 0) \, . $$
Then $a_0=1$. We prove that  $\cR'$ is  infinite by showing that
\begin{equation}\label{line-E:infinity}
a_n>2a_{n-1}    \quad (n \ge 1) \, .
\end{equation}
We use induction. As in \S\ref{ohohoh},  for $ n \ge 1$, let
$$\cU_n:=\Big\{\tau\in\cT_{\suc}(\cR'_{n-1}):\cI(\tau) \ \cap\bigcup_{P\in\cP_n}\Delta(P)\ne\emptyset\Big\}.$$
Then
$$
\cR'_n
 \ = \ \cR_{\suc}(\cR'_{n-1})\setminus\cU_n,
$$
and it follows that
\begin{equation}\label{line-E:infinity1}
a_n\ge 6a_{n-1}-\#\cU_n.
\end{equation}
On the other hand, as in \S\ref{ohohoh}, we have that
\begin{eqnarray*}
\cU_n & = & \bigcup_{k=1}^n\Big\{\tau\in\cT_{\suc}(\cR'_{n-1}):\cI(\tau) \ \cap\bigcup_{P\in\cP_{n,k}}\Delta(P)\ne\emptyset\Big\}\\
& \subset & \bigcup_{k=1}^n\bigcup_{~ \ \tau'\in\cR'_{n-k}} \!\!\! \Big\{\tau\in\cT_n:\cI(\tau) \ \cap\bigcup_{P\in\cP_{n,k}(\tau')}\Delta(P)\ne\emptyset\Big\}.
\end{eqnarray*}
Thus,  Corollary \ref{line-C:est} implies that
\begin{equation}\label{line-E:infinity2}
\#\cU_n\le \sum_{k=1}^n2a_{n-k}.
\end{equation}
On combining  \eqref{line-E:infinity1} and \eqref{line-E:infinity2}, we obtain that
\begin{equation}\label{line-E:infinity3}
a_n\ge 6a_{n-1}-\sum_{k=1}^n2a_{n-k}.
\end{equation}
With $n=1$ in \eqref{line-E:infinity3}, we  find that $a_1\ge 4$. Hence,  \eqref{line-E:infinity} holds for $n=1$. Now  assume $n\ge2$ and  that \eqref{line-E:infinity} holds with $n$  replaced by $1,\ldots,n-1$.
Then for any $1\le k\le n$, we have that  $$a_{n-k}\le 2^{-k+1}a_{n-1}.$$ Substituting this into \eqref{line-E:infinity3}, gives that
$$a_n\ge6a_{n-1}-2a_{n-1}\sum_{k=1}^n2^{-k+1}>2a_{n-1}.$$
This completes the induction step and thus establishes  \eqref{line-E:infinity}. In turn this completes the proof of Proposition \ref{P:mainlines}.
\vspace*{-2ex}

\hfill $ \boxtimes $

\vspace*{2ex}

\subsection{The inhomogeneous case: establishing Theorem \ref{thm2} \label{inhthm2strat2}}

Theorem \ref{thm2} is easily deduced from the  following  statement.

\begin{theorem}\label{T2inhomo:main}
Let  $(i,j)$ be a pair of real numbers
satisfying   $0<j\le i<1$ and $i+j=1$.   Given $a,b\in \RR$, suppose there exists
$\epsilon>0$ such that \eqref{diocondbetterT2} is satisfied.
Then $\Bad^f_{\bm\theta}(i,j)$
is a $1/2$-winning subset of $\RR$.
\end{theorem}

\noindent We have already established  the homogeneous case ($\gamma=\delta=0$) of the Theorem~\ref{T2inhomo:main};  namely Theorem \ref{T2homo:main}.  With reference to \S\ref{thm2strat2},  the crux of the `homogeneous' proof involved constructing a partition $\cP_{n } $   ($ n \ge 1$)  of $\cP$  (given by  Lemma \ref{L:partitionlines})  such that the subtree $\cS$ of an $[R]$-regular rooted tree $\cT$
has an $([R]-5)$-regular subtree $\cS'$ -- the substance of Proposition \ref{P:mainlines}.
To prove  Theorem~\ref{T2inhomo:main},  the idea is to merge the inhomogeneous constraints into the homogeneous construction as in the case of curves in \S\ref{inhThm1.2}.  More precisely, we show that $\cS'$ has an $([R]-7)$-regular subtree $\cQ'$ that incorporates the inhomogeneous constraints.   With this in mind, let
$$c':=\frac{1}{10}cR^{-2}$$
and
$$\sV:=\Big\{(p,r,q)\in\ZZ^2\times\NN : \Big|f\Big(\frac{p+\gamma}{q}\Big)-\frac{r+\delta}{q}\Big|<\frac{\kappa c'}{q^{1+j}}\Big\}.$$
Furthermore, for each $v=(p,r,q)\in\ZZ^2\times\NN$, we associate the interval $$\Delta_{\bm\theta}(v):=\Big\{x\in \RR:\Big|x-\frac{p+\gamma}{q}\Big|<\frac{c'}{q^{1+i}}\Big\}.$$
Then,  with $\cA_0 \subset \cB_0 $ as in  \S\ref{thm2strat2}, the following is the inhomogeneous analogue of Lemma~\ref{L:Bad_clines}.

\begin{lemma}\label{line-L:Bad_c'}
Let $\cA_0$, $  \sV $ and $ \Delta_{\bm\theta}(v) $ be as above.  Then
$$\cA_0\setminus\bigcup_{v\in\sV}\Delta_{\bm\theta}(v) \ \subset \ \Bad_{\bm\theta}^f(i,j)  \, . $$
\end{lemma}

\noindent{\em Proof. \ } The proof is similar to the homogeneous proof but is included for sake of completeness.
Let $x\in\cA_0$. Suppose $x\notin\Bad^f_{\bm\theta}(i,j)$. Then there exists $v=(p,r,q)\in\ZZ^2\times\NN$ such that $$\Big|x-\frac{p+\gamma}{q}\Big|<\frac{c'}{q^{1+i}}, \qquad \Big|f(x)-\frac{r+\delta}{q}\Big|<\frac{c'}{q^{1+j}}.$$ It follows that
\begin{eqnarray*}
\Big|f\Big(\frac{p+\gamma}{q}\Big)-\frac{r+\delta}{q}\Big| & \le & \Big|f\Big(\frac{p+\gamma}{q}\Big)-f(x)\Big|
+\Big|f(x)-\frac{r+\delta}{q}\Big|\\[1ex]
& < & |a|\frac{c'}{q^{1+i}}+\frac{c'}{q^{1+j}}  \ \le  \ \frac{\kappa c'}{q^{1+j}}.
\end{eqnarray*}
Thus $x\in\Delta_{\bm\theta}(v)$ and $v\in\sV$.  This completes the proof of the lemma.
\vspace*{-2ex}

\hfill $ \boxtimes $

\vspace*{2ex}

For $n\ge1$, let $$H'_n:=2c'l^{-1}R^n$$ and
$$\sV_n:=\{(p,r,q)\in\sV: H'_n\le q^{1+i}< H'_{n+1}\}.$$
Observe that  $H'_1=2c'l^{-1}R\le1$ and so it follows that $\sV=\bigcup_{n=1}^\infty\sV_n$. We inductively define a subtree $\cQ$ of $\cS'$ as follows.
Let $\cQ_0=\{\tau_0\}$. If $\cQ_{n-1}$ $(n\ge1)$ is defined, we let
$$\cQ_n:=\Big\{\tau\in\cS'_{\suc}(\cQ_{n-1}):\cI(\tau) \ \cap\bigcup_{v\in\sV_n}\Delta_{\bm\theta}(v)=\emptyset\Big\}.$$
Then $$\cQ:=\bigcup_{n=0}^\infty\cQ_n$$ is a subtree of $\cS'$ and by
construction
\begin{equation*}\label{line-IE:tau}
\cI(\tau) \ \subset \ \cA_0\setminus\bigcup_{v\in\sV_n}\Delta_{\bm\theta}(v)   \qquad \forall  \ n\ge1 \quad {\rm and } \quad  \tau\in\cQ_n.
\end{equation*}

Armed with the following result, the same arguments as in \S\ref{howzthat} with the most obvious modifications enables us to prove Theorem \ref{T2inhomo:main}. In view of this the details of the proof of Theorem \ref{T2inhomo:main} modulo Proposition \ref{line-keyinhom} are omitted.

\begin{proposition} \label{line-keyinhom}
The tree $\cQ$ has an $([R]-7)$-regular subtree.
\end{proposition}

In order to establish  the proposition, it suffices to prove the following statement.

\begin{lemma}
For any $n\ge1$ and $\tau\in\cQ_{n-1}$, there is at most one $v\in\sV_n$ such that $\cI(\tau)\cap\Delta_{\bm\theta}(v)\ne\emptyset$. Moreover, $\rho(\Delta_{\bm\theta}(v))\le lR^{-n}$. Therefore,
$$\#\Big\{\tau'\in\cS'_{\suc}(\tau):\cI(\tau') \ \cap\bigcup_{v\in\sV_n}\Delta_{\bm\theta}(v)\ne\emptyset\Big\} \ \le \ 2.$$
\end{lemma}

\noindent{\em Proof. \ }
Suppose $v_s=(p_s,r_s,q_s)\in\sV_n$ and $\cI(\tau)\cap\Delta_{\bm\theta}(v_s)\ne\emptyset$, $s=1,2$. We need to prove that $v_1=v_2$. Without loss of generality, assume that $q_1\ge q_2$. Let $x_1\in\cI(\tau)\cap\Delta_{\bm\theta}(v_1)$. The same arguments as in the proofs of \eqref{w1} and \eqref{w2} show that
\begin{equation}\label{line-w1}
|(q_1-q_2)x_1-(p_1-p_2)| \ \le \ q_2lR^{-n+1}+\frac{2c'}{q_2^i}
\end{equation}
and
\begin{equation}\label{line-w2}
\Big| (q_1-q_2)f\Big(\frac{p_1+\gamma}{q_1}\Big)- (r_1-r_2) \Big| \ \le \ q_2|a|lR^{-n+1}+\frac{4\kappa c'}{q_2^j}.
\end{equation}

Suppose for the moment  that $q_1>q_2$  and let
$$P_0:=\Big(\frac{p_1-p_2}{q_1-q_2},\frac{r_1-r_2}{q_1-q_2}\Big)  \, . $$  We show that $\cI(\tau)\cap\Delta(P_0)\ne\emptyset$ and that $P_0\in\cP$. Similar to the proof of Lemma \ref{L:inhomo-curve}, it follows from \eqref{line-w1} that
$$(q_1-q_2)^{1+i}\Big|x_1-\frac{p_1-p_2}{q_1-q_2}\Big| \ <  \ c.$$
So $x_1\in\Delta(P_0)$ and it follows that  $\cI(\tau)\cap\Delta(P_0)\ne\emptyset$. Moreover, by making use
of \eqref{line-w2} we have that
$$(q_1-q_2)^{1+j} \Big|f\Big(\frac{p_1-p_2}{q_1-q_2}\Big)-\frac{r_1-r_2}{q_1-q_2}\Big| \ <  \ \kappa c.$$
Thus $P_0\in\cP$ and so there exists a  unique integer $n_0\ge1$  such that $P_0\in\cP_{n_0}$.   The same argument as in the proof of Lemma \ref{L:inhomo-curve} shows that $n_0\ge n $, and so
$$(q_1-q_2)|E_{P_0}| \ \ge \  H_{n_0} \ \ge \  H_{n} \ = \  4\kappa cl^{-1}R^n. $$
On the other hand, we have that
\begin{eqnarray*}
(q_1-q_2)|E_{P_0}| & \le  &  \kappa(q_1-q_2)^{1+i} \ \le \ \kappa q_1^{1+i} \\[1ex]
& \le  &  \kappa H'_{n+1} \ = \ \frac{1}{5}\kappa cl^{-1}R^{n-1}  \, .
\end{eqnarray*}
This contradicts the above lower bound for $(q_1-q_2)|E_{P_0}|$ and we conclude that $q_1=q_2$.  It now  follows from \eqref{line-w1} and \eqref{line-w2} that
$$|p_1-p_2|<1 \qquad {\rm and }  \qquad |r_1-r_2|<1.$$
Thus,  $p_1=p_2$ and $r_1=r_2$. In other words,  $v_1=v_2$ and this proves the  main  substance  of the proposition.  The proofs of the remaining parts are the same as those for Lemma \ref{L:inhomo-curve}.
\vspace*{-1ex}

\hfill $ \boxtimes $

\vspace*{2ex}

As already mentioned, given Proposition \ref{line-keyinhom},   the proof of Theorem \ref{T2inhomo:main} follows on adapting the arguments of \S\ref{howzthat}.

\section{The proof of  Theorem \ref{2DT:main} \label{AN}}

We need the notion of  regular colorings of rooted trees.  Let $D\in\NN$. A \emph{$D$-coloring} of a rooted tree $\cT$ is a map $\gamma:\cT\to\{1,\ldots,D\}$. For $\cV\subset\cT$ and $1\le i\le D$, we denote $\cV^{(i)}=\cV\cap\gamma^{-1}(i)$. Let $N\in\NN$ be an integer multiple of $D$, and suppose that $\cT$ is $N$-regular. We say that a $D$-coloring of $\cT$ is \emph{regular} if for any $\tau\in\cT$ and $1\le i\le D$, we have $\#\cT_{\suc}(\tau)^{(i)}=N/D$. The following two types of subtrees are of interest to us.
\begin{itemize}
  \item The subtree $\cS$ is of \emph{type (I)} if for any $\tau\in\cS$ and $1\le i\le D$, we have $\#\cS_{\suc}(\tau)^{(i)}=1$.
  \item The subtree $\cS$ is of \emph{type (II)} if for any $\tau\in\cS$, there exists $1\le i(\tau)\le D$ such that $\cS_{\suc}(\tau)=\cT_{\suc}(\tau)^{(i(\tau))}$.
\end{itemize}
Roughly speaking, in the proof of Theorem \ref{2DT:main}, the two types of subtrees correspond to strategies of the two players in Schmidt's game.
We will make use of the following criterion for the existence of subtree of type (I).  It appears as Proposition 2.2 in \cite{Jinpeng2}.

\begin{proposition}\label{2DP:tree}
Let $\cT$ be an $N$-regular rooted tree with a regular $D$-coloring, and let $\cS\subset\cT$ be a subtree. Suppose that for every subtree $\cR\subset\cT$ of type (II), $\cS\cap\cR$ is infinite. Then $\cS$ contains a subtree of type (I).
\end{proposition}

\subsection{The winning strategy for Theorem \ref{2DT:main} \label{winstratinhom}}

Let $\alpha_0:=(30\sqrt{2})^{-1}$, $\beta\in(0,1)$. We want to prove that $\Bad_{\bm\theta}(i,j)$ is $(\alpha_0,\beta)$-winning. In the first round of the game, \textbf{B}hupen chooses a closed disc $\cB_0\subset \RR^2$.  Now  \textbf{A}yesha chooses any closed disc $\cA_0\subset\cB_0$ with diameter $\rho(\cA_0)=\alpha_0\rho(\cB_0)$. Let
$$l:=\rho(\cA_0), \qquad R:=(\alpha_0\beta)^{-1}, \qquad m:=15.$$
By a \emph{square} we mean a set of the form
$$\Sigma=\{(x,y)\in\RR^2:x_0\le x\le x_0+\ell(\Sigma), y_0\le y\le y_0+\ell(\Sigma)\},$$ where $\ell(\Sigma)>0$ is the side length of $\Sigma$.
Let $\Sigma_0$ be the circumscribed square of $\cA_0$. Then $\ell(\Sigma_0)=l$. Let $\cT$ be an $m^2[R/m]^2$-regular rooted tree with a regular $[R/m]^2$-coloring. We choose and fix
an injective map $\Phi$ from $\cT$ to the set of subsquares of $\Sigma_0$ satisfying the following conditions:
\begin{itemize}
  \item For any $n\ge0$ and $\tau\in\cT_n$, we have $\ell(\Phi(\tau))=lR^{-n}$. In particular, the root $\tau_0$ of $\cT$ is mapped to $\Sigma_0$.
  \item For $\tau,\tau'\in\cT$, if $\tau\prec\tau'$, then $\Phi(\tau)\subset\Phi(\tau')$.
  \item For any $n\ge1$ and $\tau\in\cT_{n-1}$, the interiors of the squares $\{\Phi(\tau'):\tau'\in\cT_{\suc}(\tau)\}$ are mutually disjoint, the union $\bigcup_{\tau'\in\cT_{\suc}(\tau)}\Phi(\tau')$ is a square of side length $m[R/m]lR^{-n}$, and for any $1\le i\le [R/m]^2$, the union $\bigcup_{\tau'\in\cT_{\suc}(\tau)^{(i)}}\Phi(\tau')$ is a square of side length $mlR^{-n}$.
\end{itemize}

Let $c>0$ be such that
\begin{equation}\label{2DE:c}
c<\min\Big\{\frac{1}{6}lR^{-1},\frac{1}{16}R^{-12}\Big\},
\end{equation}
and in turn let
\begin{equation}\label{2DE:c'}
c':=\frac{1}{6}cR^{-2}.
\end{equation}
For each  $P=(\frac{p}{q},\frac{r}{q})\in\QQ^2$, we associate the rectangle $$\Delta(P):=\Big\{(x,y)\in\RR^2:\Big|x-\frac{p}{q}\Big|\le\frac{c}{q^{1+i}},\Big|y-\frac{r}{q}\Big|\le\frac{c}{q^{1+j}}\Big\},$$
and for $v=(p,r,q)\in\ZZ^2\times\NN$, we associate the rectangle
$$\Delta_{\bm\theta}(v):=\Big\{(x,y)\in\RR^2:\Big|x-\frac{p+\gamma}{q}\Big|\le\frac{c'}{q^{1+i}},\Big|y-\frac{r+\delta}{q}\Big|\le\frac{c'}{q^{1+j}}\Big\}.$$
Then
\begin{equation}\label{2DE:Badc2}
\RR^2\setminus\Big(\bigcup_{P\in\QQ^2}\Delta(P) \ \cup\bigcup_{v\in\ZZ^2\times\NN}\Delta_{\bm\theta}(v)\Big) \ \subset \ \Bad(i,j)\cap\Bad_{\bm\theta}(i,j).
\end{equation}
For $n\ge1$, let
$$H_n:=6cl^{-1}R^n \, , \qquad H'_n:=3c'l^{-1}R^n $$ and define
\begin{equation}\label{2DE:C_n}
\cP_n:=\Big\{P=\Big(\frac{p}{q},\frac{r}{q}\Big)\in\QQ^2: H_n\le q\max\{|A_P|,|B_P|\}< H_{n+1}\Big\}
\end{equation}
and
\begin{equation}\label{2DE:C_n'}
\sV_n:=\{v=(p,r,q)\in\ZZ^2\times\NN: H'_n\le q^{1+\max\{i,j\}}<H'_{n+1}\},
\end{equation}
where $A_P$ and $B_P$ are as in \S\ref{premain}. In view of  \eqref{2DE:c}, we have that  $H'_1\le H_1=6cl^{-1}R\le1$. Thus $$\QQ^2=\bigcup_{n=1}^\infty\cP_n \quad {\rm and} \quad  \ZZ^2\times\NN=\bigcup_{n=1}^\infty\sV_n  \, . $$
We inductively define a subtree $\cS$ of $\cT$ as follows.
Let $\cS_0=\{\tau_0\}$. If $n\ge1$ and $\cS_{n-1}$ is defined, we let
\begin{equation}\label{2DE:S_n}
\cS_n:=\Big\{\tau\in\cT_{\suc}(\cS_{n-1}):\Phi(\tau) \ \cap \ \Big(\bigcup_{P\in\cP_n}\Delta(P)\cup\bigcup_{v\in\sV_n}\Delta_{\bm\theta}(v)\Big)=\emptyset\Big\}.
\end{equation}
Then $\cS=\bigcup_{n=0}^\infty\cS_n$ is a subtree of $\cT$ and by construction
\begin{equation}\label{2DE:tau}
\Phi(\tau) \ \subset \ \RR^2\setminus\Big(\bigcup_{P\in\cP_n}\Delta(P)\cup\bigcup_{v\in\sV_n}\Delta_{\bm\theta}(v)\Big) \qquad \forall \  n\ge1 \quad {\rm and} \quad  \tau\in\cS_n.
\end{equation}

\medskip

The following proposition is the key to proving Theorem \ref{2DT:main}.

\begin{proposition}\label{2DP:main}
The tree $\cS$ contains a subtree of type (I).
\end{proposition}

\subsubsection{Proof of Theorem \ref{2DT:main} modulo Proposition \ref{2DP:main}}

This is essentially the same as the proof of Theorem 1.1 in \cite{Jinpeng2}.  However for the sake of completeness we have included the short argument.
Let $\cS'$ be a subtree of $\cS$ of type (I). We inductively prove that for every $n\ge0$,
\begin{equation}\label{2DE:state}
\text{\textbf{A}yesha  can choose $\cA_{n}$ to be the inscribed closed disc of $\Phi(\tau_n)$ for some $\tau_n\in\cS'_n$.}
\end{equation}
If $n=0$, there is nothing to prove. Assume $n\ge1$ and \textbf{A}yesha  has chosen $\cA_{n-1}$ as the inscribed closed disc of $\Phi(\tau_{n-1})$, where $\tau_{n-1}\in\cS'_{n-1}$. For any choice $\cB_n\subset\cA_{n-1}$ of \textbf{B}hupen, the inscribed square of $\cB_n$ has side length
$$\frac{\sqrt{2}}{2}\rho(\cB_n) \ = \ \frac{\sqrt{2}}{2}\beta\rho(\cA_{n-1}) \ = \ \frac{\sqrt{2}}{2}\beta\ell(\Phi(\tau_{n-1})) \ = \ \frac{\sqrt{2}}{2}\beta lR^{-n+1} \ = \ 2mlR^{-n}.$$
Thus there exists $1\le i\le [R/m]^2$ such that $\bigcup_{\tau\in\cT_{\suc}(\tau_{n-1})^{(i)}}\Phi(\tau)\subset\cB_n$. Let $\tau_n$ be the unique vertex in $\cS'_{\suc}(\tau_{n-1})^{(i)}$. Then $\Phi(\tau_n)\subset\cB_n$. Note that the diameter of the inscribed closed disc of $\Phi(\tau_n)$ is equal to
$$\ell(\Phi(\tau_n)) \ = \ R^{-1}\ell(\Phi(\tau_{n-1})) \ = \ \alpha_0\beta\rho(\cA_{n-1}) \ = \ \alpha_0\rho(\cB_n).$$
So \textbf{A}yesha  can choose $\cA_n$ to be the inscribed closed disc of $\Phi(\tau_n)$. This proves \eqref{2DE:state}.

In view of \eqref{2DE:state}, \eqref{2DE:tau} and \eqref{2DE:Badc2}, we have
\begin{eqnarray*}
\bigcap_{n=0}^\infty \cA_n \ \subset \ \bigcap_{n=0}^\infty\Phi(\tau_n) & \subset & \bigcap_{n=1}^\infty\RR^2\setminus\Big(\bigcup_{P\in\cP_n}\Delta(P)\cup\bigcup_{v\in\sV_n}\Delta_{\bm\theta}(v)\Big)\\
& = & \RR^2\setminus\Big(\bigcup_{P\in\QQ^2}\Delta(P) \ \cup\bigcup_{v\in\ZZ^2\times\NN}\Delta_{\bm\theta}(v)\Big) \ \subset \ \Bad(i,j)\cap\Bad_{\bm\theta}(i,j).
\end{eqnarray*}
This proves the theorem assuming the truth of Proposition \ref{2DP:main}.
\hfill $ \boxtimes $

\subsection{Proof of  Proposition \ref{2DP:main}}
Let $w>0$. By a \emph{strip} of width $w$, we mean a subset of $\RR^2$ of the form
$$\mathcal{L}:=\{\mathbf{x}\in\RR^2:|\mathbf{x}\cdot\mathbf{u}-a|\le w/2\},$$
where the dot denotes the standard inner product, $\mathbf{u}\in\RR^2$ is a unit vector, and $a\in\RR$.
The following result is proved in \cite[Corollary 4.2]{Jinpeng2}.

\begin{lemma}\label{2DL:strip}
For any $n\ge1$, there exists a partition $\cP_n=\bigcup_{k=1}^n\cP_{n,k}$ such that for any $1\le k\le n$ and $\tau\in\cS_{n-k}$, there is a strip of width $\frac{2}{3}lR^{-n}$ which contains all the rectangles $$\{\Delta(P):P\in\cP_{n,k},\Phi(\tau)\cap\Delta(P)\ne\emptyset\}.$$
\end{lemma}

We now prove a corresponding result which takes into consideration the  inhomogeneous approximation aspect.

\begin{lemma}\label{2DL:strip2}
For any $n\ge1$ and $\tau\in\cS_{n-1}$, there is at most one $v\in\sV_n$ such that $\Phi(\tau)\cap\Delta_{\bm\theta}(v)\ne\emptyset$. Moreover, $\Delta_{\bm\theta}(v)$ is contained in a strip of width $\frac{2}{3}lR^{-n}$.
\end{lemma}

\noindent {\em Proof. \ }
Suppose $v_s=(p_s,r_s,q_s)\in\sV_n$ and $\Phi(\tau)\cap\Delta_{\bm\theta}(v_s)\ne\emptyset$, $s=1,2$. We need to prove that $v_1=v_2$. Without loss of generality, assume that $q_1\ge q_2$. Since $\Phi(\tau)\cap\Delta_{\bm\theta}(v_s)\ne\emptyset$, there exists $(x_s,y_s)\in\Phi(\tau)$ such that
$$|q_sx_s-(p_s+\gamma)|\le\frac{c'}{q_s^i} \, , \qquad |q_sy_s-(r_s+\delta)|\le\frac{c'}{q_s^j}.$$
It follows that
\begin{eqnarray}
|(q_1-q_2)x_1-(p_1-p_2)| & \le & q_2|x_1-x_2|+|q_1x_1-(p_1+\gamma)|+|q_2x_2-(p_2+\gamma)| \nonumber \\[1ex]
& \le & q_2lR^{-n+1}+\frac{2c'}{q_2^i}. \label{yu}
\end{eqnarray}
Similarly, we have that
\begin{equation} \label{yu2}
|(q_1-q_2)y_1-(r_1-r_2)| \ \le \ q_2lR^{-n+1}+\frac{2c'}{q_2^j}.
\end{equation}

We first prove that $q_1=q_2$. Suppose this is not the case. Then
\begin{eqnarray*}
\Big|x_1-\frac{p_1-p_2}{q_1-q_2}\Big| & \le & \frac{1}{q_1-q_2}\Big(q_2lR^{-n+1}+\frac{2c'}{q_2^i}\Big)\\
 & \le & \frac{1}{(q_1-q_2)^{1+i}}\Big(q_1^iq_2lR^{-n+1}+2c'\frac{q_1^i}{q_2^i}\Big)\\
 & \le & \frac{1}{(q_1-q_2)^{1+i}}(H'_{n+1}lR^{-n+1}+2c'R)\\
& = & \frac{1}{(q_1-q_2)^{1+i}}\Big(\frac{c}{2}+\frac{1}{3}cR^{-1}\Big)\\
 & \le & \frac{c}{(q_1-q_2)^{1+i}}.
\end{eqnarray*}
Similarly, $$\Big|y_1-\frac{r_1-r_2}{q_1-q_2}\Big| \ \le \ \frac{c}{(q_1-q_2)^{1+j}}.$$
Thus,  if we let $$P_0:=\Big(\frac{p_1-p_2}{q_1-q_2},\frac{r_1-r_2}{q_1-q_2}\Big),$$ then $(x_1,y_1)\in\Delta(P_0)$. In particular, $\Phi(\tau)\cap\Delta(P_0)\ne\emptyset$.
Let $n_0\ge1$ be the unique integer such that $P_0\in\cP_{n_0}$. It is easily verified that
\begin{equation}\label{2DE:m}
n_0\ge n.
\end{equation}
Indeed, if $n_0\le n-1$, then $\cS_{n_0}$ contains an ancestor $\tau'$ of $\tau$ and by \eqref{2DE:S_n} we have  that
$$\Phi(\tau)\cap\Delta(P_0) \, \subset  \, \Phi(\tau')\cap\Delta(P_0)=\emptyset \, .$$
This is a  contradiction since the left hand side is non-empty.
Now, in view of \eqref{2DE:C_n} and \eqref{2DE:m}, we have that
$$(q_1-q_2)^{1+\max\{i,j\}} \ \ge \  (q_1-q_2)\max\{|A_{P_0}|,|B_{P_0}|\} \ \ge \  H_{n_0} \ \ge \  H_n \ = \ 6cl^{-1}R^n.$$
On the other hand, we have that
$$(q_1-q_2)^{1+\max\{i,j\}} \ \le \  q_1^{1+\max\{i,j\}} \ \le \  H'_{n+1} \ = \ \frac{1}{2}cl^{-1}R^{n-1} \, .$$
This contradicts the above lower bound and so we must have that $q_1=q_2$. It then follows via (\ref{yu}) and (\ref{yu2}) that
$$|p_1-p_2| \ \le \  q_2lR^{-n+1}+2c' \ \le \  c+2c' \ < \ 1,$$
$$|r_1-r_2| \ \le \  q_2lR^{-n+1}+2c' \ \le \  c+2c' \ < \ 1.$$
The left hand sides of these inequalities are integers so we must have that $p_1=p_2$ and $r_1=r_2$. The upshot of this is that  $v_1=v_2$ and  so establishes  the main substance of the lemma.  Regarding the moreover part, simply observe that  for any $v=(p,r,q)\in\sV_n$, $\Delta_{\bm\theta}(v)$ is contained in a strip of width
$$\min\Big\{\frac{2c'}{q^{1+i}},\frac{2c'}{q^{1+j}}\Big\} \ = \ \frac{2c'}{q^{1+\max\{i,j\}}} \ \le \  \frac{2c'}{H'_n} \ = \ \frac{2}{3}lR^{-n}.$$
\hfill $ \boxtimes $

\vspace{2ex}

The following result proved in  \cite[Lemma 4.3]{Jinpeng2} gives an upper bound for the number of certain squares which
intersect a thin strip.

\begin{lemma}\label{2DL:sub}
Let $\cR\subset\cT$ be a subtree of type (II), let $n\ge1$, and let $\mathcal{L}$ be a strip of width $\frac{2}{3}lR^{-n}$. Then for any $1\le k\le n$ and $\tau\in\cR_{n-k}$, we have that
$$\#\{\tau'\in\cR(\tau)_k:\Phi(\tau')\cap\mathcal{L}\ne\emptyset\} \ \le \ (3m-2)^k.$$
\end{lemma}

On combining Lemmas \ref{2DL:strip}, \ref{2DL:strip2} and \ref{2DL:sub}, we obtain the following statement.

\begin{corollary}\label{2DC:<=2}
Let $\cR\subset\cT$ be a subtree of type (II) and let $n\ge1$. Then
\begin{itemize}
  \item For any $\tau\in\cS_{n-1}\cap\cR_{n-1}$, we have
  $$\#\Big\{\tau'\in\cR_{\suc}(\tau):\Phi(\tau') \ \cap \ \Big(\bigcup_{P\in\cP_{n,1}}\Delta(P)\cup\bigcup_{v\in\sV_n}\Delta_{\bm\theta}(v) \, \Big)\ne\emptyset\Big\} \ \le \ 2(3m-2).$$
  \item For any $2\le k\le n$ and $\tau\in\cS_{n-k}\cap\cR_{n-k}$, we have
$$\#\Big\{\tau'\in\cR(\tau)_k:\Phi(\tau') \ \cap\bigcup_{P\in\cP_{n,k}}\Delta(P)\ne\emptyset\Big\} \ \le \ (3m-2)^k.$$
\end{itemize}
\end{corollary}

We are now in the position to prove Proposition \ref{2DP:main}. The proof  is essentially the same as the proof of Proposition 3.3  in \cite{Jinpeng2}.  However for the sake of completeness we have included the argument.
In view of Proposition \ref{2DP:tree}, it suffices to prove that the intersection of $\cS$ with every subtree of type (II) is infinite. Let $\cR\subset\cT$ be a subtree of type (II), and let
$$\cR':=\cR\cap\cS,    \qquad {\rm   and }  \qquad  a_n:=\#\cR'_n \qquad (n \ge 0 ) \, . $$
Then $a_0=1$. We prove that $\cR'$ is infinite  by showing that
\begin{equation}\label{2DE:infinity}
a_n>112a_{n-1}   \quad (n \ge 1) \, .
\end{equation}
We use induction.  For $n \ge 1$, let
$$\cU_n:=\Big\{\tau\in\cR_{\suc}(\cR'_{n-1}):\Phi(\tau) \ \cap \ \Big(\bigcup_{P\in\cP_n}\Delta(P)\cup\bigcup_{v\in\sV_n}\Delta_{\bm\theta}(v)\Big)\ne\emptyset\Big\}.$$
It is easy to see from \eqref{2DE:S_n} that
$$
\cR'_n = \cR_{\suc}(\cR'_{n-1}) \setminus \cU_n \, ,
$$
and so it follows  that
\begin{equation}\label{2DE:infinity1}
a_n  \, = \, \#\cR_{\suc}(\cR'_{n-1})-\#\cU_n  \, =  \, m^2a_{n-1}-\#\cU_n.
\end{equation}
On the other hand, we have that
\begin{eqnarray*}
\cU_n & = & \Big\{\tau'\in\cR_{\suc}(\cR'_{n-1}):\Phi(\tau') \ \cap \ \Big(\bigcup_{P\in\cP_{n,1}}\Delta(P)\cup\bigcup_{v\in\sV_n}\Delta_{\bm\theta}(v)\Big)\ne\emptyset\Big\}\\
& & \qquad  \cup   \quad \bigcup_{k=2}^n
\Big\{\tau'\in\cR_{\suc}(\cR'_{n-1}):\Phi(\tau') \ \cap\bigcup_{P\in\cP_{n,k}}\Delta(P)\ne\emptyset\Big\}\\[2ex]
& \subset & \bigcup_{\tau\in\cR'_{n-1}}
\Big\{\tau'\in\cR_{\suc}(\tau):\Phi(\tau') \ \cap \ \Big(\bigcup_{P\in\cP_{n,1}}\Delta(P)\cup\bigcup_{v\in\sV_n}\Delta_{\bm\theta}(v)\Big)\ne\emptyset\Big\}\\
& &  \qquad \cup   \quad \bigcup_{k=2}^n   \ \ \bigcup_{\tau\in\cR'_{n-k}}
\Big\{\tau'\in\cR(\tau)_k:\Phi(\tau') \ \cap\bigcup_{P\in\cP_{n,k}}\Delta(P)\ne\emptyset\Big\}.
\end{eqnarray*}
Thus, Corollary \ref{2DC:<=2} implies that
\begin{equation}\label{2DE:infinity2}
\#\cU_n \ \le \ 2(3m-2)a_{n-1}+\sum_{k=2}^n(3m-2)^ka_{n-k}.
\end{equation}
On combining \eqref{2DE:infinity1} and \eqref{2DE:infinity2}, we obtain that
\begin{equation}\label{2DE:infinity3}
a_n \ \ge \  (m^2-3m+2)a_{n-1}-\sum_{k=1}^n(3m-2)^ka_{n-k} \ = \ 182a_{n-1}-\sum_{k=1}^n43^ka_{n-k}.
\end{equation}
With $n=1$ in \eqref{2DE:infinity3}, we find that $a_1\ge139$. Hence,  \eqref{2DE:infinity} holds for $n=1$. Now assume $n\ge2$ and that  \eqref{2DE:infinity} holds  with  $n$ replaced by $1,\ldots,n-1$.
Then for any $1\le k\le n$, we have that  $$a_{n-k}\le112^{-k+1}a_{n-1}.$$ Substituting this into \eqref{2DE:infinity3}, gives that
$$a_n \ \ge \  \Big(182-112\sum_{k=1}^n(43/112)^k\Big)a_{n-1} \ > \ 112a_{n-1}.$$
This completes the induction step and thus establishes \eqref{2DE:infinity}. In turn this completes the proof of Proposition \ref{2DP:main}.
\hfill $ \boxtimes $

\bigskip
\bigskip

\noindent\textit{Acknowledgements.}
This paper was (finally) completed during the authors' stay at the Isaac Newton Institute for Mathematical Sciences at Cambridge during the programme ``Interactions between Dynamics of Group Actions and Number Theory'' 9 June -- 4 July 2014. We thank the organisers for  inviting us  and for creating a wonderful and relaxed working environment. JA would like to thank the great hospitality provided by the University of York where this work was initiated during his stay in June -- July of 2012. SV would like to use this opportunity to acknowledge his school friend and fellow mathematician Sunil Talwar (1964-2014) who died very unexpectedly just short of his fiftieth birthday.  We had many memorable times together and shared many lifetimes of laughter  -- even when discussing analysis.  Thank you Sunil!

\vspace{15mm}

\vspace{15mm}

\noindent Jinpeng An: LMAM, School of Mathematical Sciences,

\vspace{0mm}

\noindent\phantom{Jinpeng An : }Peking University, Beijing,  100871, China.


\noindent\phantom{Jinpeng An : }e-mail: anjinpeng@gmail.com

\vspace{3mm}

\noindent Victor Beresnevich: Department of Mathematics, University of
York,

\vspace{0mm}

\noindent\phantom{Victor Beresnevich : }Heslington, York, YO10 5DD,
England.


\noindent\phantom{Victor Beresnevich : }e-mail: vb8@york.ac.uk

\vspace{3mm}

\noindent Sanju Velani: Department of Mathematics, University of
York,

\vspace{0mm}

\noindent\phantom{Sanju Velani : }Heslington, York, YO10 5DD,
England.


\noindent\phantom{Sanju  Velani : }e-mail: slv3@york.ac.uk

\end{document}